\documentclass[11pt]{article}

\usepackage[latin1]{inputenc}
\usepackage[english]{babel}
\usepackage{amsfonts,amssymb,amsmath}
\usepackage{mathrsfs}
\usepackage{mathtext}
\usepackage{textcomp}
\usepackage{graphicx}
\usepackage{epsfig}
\usepackage{rotating}

\usepackage{tikz}
\usetikzlibrary{matrix,arrows,decorations.pathmorphing,backgrounds,fit}
\usepackage{subfig}

\usepackage[toc,page]{appendix}

\usepackage[all]{xy}

\usepackage{hyperref}

\newcommand{\lra}{\longrightarrow}
\newcommand{\sm}{\setminus}

\newcommand{\C}{\mathbb{C}}
\newcommand{\R}{\mathbb{R}}
\newcommand{\Z}{\mathbb{Z}}

\renewcommand{\H}{\mathcal{H}}
\newcommand{\Hyp}{\mathbb{H}}

\newcommand{\J}{\mathcal{J}}
\newcommand{\Hb}{\mathfrak{H}}

\newcommand{\CP}{\mathbb{CP}}
\newcommand{\De}{\mathcal{D}}
\renewcommand{\P}{\mathscr{P}}

\newcommand{\Sp}{\mathrm{Sp}}
\newcommand{\SL}{\mathrm{SL}}

\newcommand{\di}{\mathbf{d}}
\newcommand{\Id}{\mathrm{Id}}
\newcommand{\inter}{\mathrm{int}}

\newcommand{\card}{\mathrm{Card}}
\newcommand{\lbd}{\lambda}

\newcommand{\dem}{\noindent {\textbf{Proof:} }}
\newcommand{\rem}{ \noindent \textbf{Remark: }}
\newcommand{\ex}{ \noindent \textbf{Example: }}

\renewcommand{\geq}{\geqslant}
\renewcommand{\leq}{\leqslant}
\newcommand{\transpose}[1]{{\vphantom{#1}}^{\mathit t}{#1}}

\newcommand{\vide}{\varnothing}
\newcommand{\carre}{\hfill $\Box$\\}

\newcommand{\ra}{\rightarrow}
\newcommand{\DS}{\displaystyle}

\newtheorem{theorem}{Theorem}[section]
\newtheorem{definition}[theorem]{Definition}
\newtheorem{proposition}[theorem]{Proposition}
\newtheorem{prop-def}[theorem]{Proposition-definition}
\newtheorem{corollary}[theorem]{Corollary}
\newtheorem{lemma}[theorem]{Lemma}

\begin{document}

\title{\bf On the topology of $\mathcal{H}(2)$}


\author{ Duc-Manh NGUYEN\\ Institut de Math\'ematiques de Bordeaux\\ Universit\'e Bordeaux 1\\ 351, cours de la Lib\'eration\\ 33405 Talence Cedex\\ France \\ \texttt{duc-manh.nguyen@math.u-bordeaux1.fr}}


\date{March 8, 2012}

\maketitle

\begin{abstract}

The space $\H(2)$ consists of pairs $(M,\omega)$, where $M$ is a Riemann surface of genus two, and $\omega$ is a holomorphic $1$-form which has
only one zero of order two. There exists a natural action of $\C^*$ on $\H(2)$ by multiplication to the holomorphic $1$-form. In this paper, we
single out a proper subgroup $\Gamma$ of $\Sp(4,\Z)$ generated by three elements, and show that the space $\H(2)/\C^*$ can be identified with the
quotient $\Gamma\backslash\J_2$, where $\J_2$ is the Jacobian locus in the Siegel upper half space $\mathfrak{H}_2$. A direct consequence of this result is that $[\mathrm{Sp}(4,\Z):\Gamma]=6$. The group $\Gamma$ can also be interpreted as the image of the fundamental group of $\H(2)/\C^*$ in the symplectic group $\Sp(4,\Z)$. 

\end{abstract}

\section{Introduction}

In this paper we are concerned with translation surfaces in the stratum $\H(2)$. Each element of $\H(2)$ can be either considered as a
translation surface having only one singularity of angle $6\pi$ together with a parallel line field, or a pair $(M,\omega)$, where $M$ is a
Riemann surface of genus two, and $\omega$ is a holomorphic $1$-form having a single zero of order two on $M$. Using the latter viewpoint, we
see that $\C^*$ acts naturally on $\H(2)$ by multiplication to the holomorphic $1$-form. Note that, if $\omega$ has only one zero on $M$, then
this zero must be a Weierstrass point of $M$. Therefore, the quotient $\H(2)/\C^*$ consists of pairs $(M,W)$, where $M$ is a Riemann surface of
genus two, and $W$ is a Weierstrass point of $M$, two pairs $(M_1,W_1)$ and $(M_2,W_2)$ are identified if there exists a conformal homeomorphism
$\phi: M_1 \lra M_2$ such that $\phi(W_1)=W_2$.

\noindent The space $\H(2)$ is well known to be a complex orbifold of dimension $4$. For any pair $(M,\omega)$, let $\gamma_1,\dots,\gamma_4$ be
a basis of the group $H_1(M,\Z)$, then the {\em period mapping}

$$(M,\omega) \longmapsto (\int_{\gamma_1}\omega,\dots,\int_{\gamma_4}\omega) \in \C^4$$

\noindent  gives a local chart for $\H(2)$ in a neighborhood of $(M,\omega)$. Consequently, we see that $\H(2)/\C^*$ can be  endowed with a
complex projective orbifold structure.

In Section \ref{GrpGamSect}, we introduce a construction of translations surface in $\H(2)$ from triples of parallelograms by a unique gluing model. This construction gives rise to the notion of {\em parallelogram decomposition} of surfaces in $\H(2)$. Actually, given a translation
surface $(M,\omega)$ in $\H(2)$, there exist infinitely many parallelogram decompositions of $(M,\omega)$. From a fixed parallelogram
decomposition, one can obtain other ones by applying some elementary moves,  which are called $T,S$, and $R$, these moves are realized by some
homeomorphisms of the surface $M$ whose actions on the group $H_1(M,\Z)$ (in an appropriate basis) are given by the following matrices
respectively

$$T=\left(%
\begin{array}{cccc}
  1 & 1 & 0 & 0 \\
  0 & 1 & 0 & 0 \\
  0 & 0 & 1 & 0 \\
  0 & 0 & 0 & 1 \\
\end{array}%
\right), \; S=\left(%
\begin{array}{cccc}
  0 & 1 & 0 & 1 \\
  0 & 0 & -1 & 0 \\
  0 & 1 & 0 & 0 \\
  -1 & 0 & 1 & 0 \\
\end{array}%
\right), \; R=\left(%
\begin{array}{cccc}
  1 & 0 & 0 & 0 \\
  0 & 1 & 0 & 0 \\
  0 & 0 & 1 & 1 \\
  0 & 0 & 0 & 1 \\
\end{array}%
\right).$$

Let $\Gamma$ be the group generated by the matrices $T,S$, and $R$, then $\Gamma$ is a proper subgroup of $\mathrm{Sp}(4,\Z)$. The main result
of this paper is the following:

\begin{theorem}\label{th0}

There exists a homeomorphism from $\H(2)/\C^*$ to the quotient $\Gamma\backslash\J_2$, where $\J_2$ is the Jacobian locus of Riemann surfaces of
genus two in the Siegel upper half space $\Hb_2$. As a consequence, we have $[\mathrm{Sp}(4,\Z):\Gamma]=6$.

\end{theorem}

Since every Riemann surface of genus two is a two-sheeted branched cover of the sphere $\CP^1$, there exists a natural identification between $\H(2)/\C^*$ and the moduli space $S^*_{0,5}$ of pairs $(\lbd_0,\underline{\lbd})$, where $\lbd_0\in \CP^1$, and $\underline{\lbd}$ is a subset of cardinal five of $\CP^1\sm \{\lbd_0\}$.

Let $\rm{Mod}_{0,6}$ denote the mapping class group of the sphere with six punctures, and $\rm{Mod}^*_{0,5}$, the subgroup of index $6$ in $\rm{Mod}_{0,6}$ that fixes one of the punctures. The space $\H(2)/\C^*$ can be identified with the quotient $\mathcal{T}_{0,6}/\rm{Mod}_{0,5}^*$, where $\mathcal{T}_{0,6}$ is the Teichm\"{u}ller space of the sphere with six punctures. The group $\Gamma$ can be then considered as the image of $\rm{Mod}^*_{0,5}$ in $\Sp(4,\Z)$. Note that we have an isomorphism between $\rm{Mod}^*_{0,5}$ and $B_5/Z(B_5)$, where $B_5$ is the braid group of the closed disc with five punctures, and $Z(B_5)$ is the center of $B_5$.

It is well known that the Jacobian locus $\mathcal{J}_2$ is the quotient of $\mathcal{T}_{0,6}$ by the Torelli group $\mathcal{I}_2$. Thus we have the following diagram
\begin{center}
\begin{tikzpicture}[scale=1]
\matrix(m)[matrix of math nodes, row sep=2em, column sep=2em]{ \mathcal{T}_{0,6} &  \\  & \mathcal{J}_2 \cong \mathcal{T}_{0,6}/\mathcal{I}_2 \\
 \H(2)/\C^*\cong \mathcal{T}_{0,6}/{\rm Mod}_{0,5}^* & \\ };
\path[->, >=angle 90] (m-1-1) edge (m-2-2);
\path[->, >=angle 90] (m-1-1) edge (m-3-1);
\path[->, >=angle 90, dashed] (m-2-2) edge (m-3-1);

\end{tikzpicture}
\end{center}

Since the action of $-\Id$ on $\mathfrak{H}_2$ is trivial, a direction consequence of Theorem~\ref{th0} is the following

\begin{corollary}\label{cor0}
We have the following exact sequence

\begin{equation}\label{MainExtSeq}
1 \rightarrow \mathcal{I}_2 \rightarrow {\rm Mod}^*_{0,5} \cong  B_5/Z(B_5) \rightarrow \Gamma/\{\pm \Id\} \rightarrow 1
\end{equation} 

\end{corollary}

\bigskip

We refer to \cite{McM09} for a more detailed account on symplectic and unitary representations of braid groups.

\bigskip

The paper is organized as follows: in Section~\ref{ThetSect}, we recall the basic properties of the Theta functions, and give a brief explanation of how these functions allow us to compute the branched points of hyperelliptic coverings. In the following section, we introduce the notion of {\em parallelogram decomposition} for surfaces in $\H(2)$, and the three elementary moves on those decompositions. We then define the group $\Gamma$ as the group generated by the homology action of the elementary moves. Note that there are parallelogram decompositions for which some elementary moves can not be carried out. This means that the set of parallelogram decompositions is not the right place to study the action of the group generated by the elementary moves. To fix this problem, in Section~\ref{AdmDecSect}, we introduce the notion of {\em admissible decomposition}, which generalizes the one of parallelogram decomposition. For every admissible decomposition, all the elementary moves can be carried out. The key ingredient of the proof of Theorem~\ref{th0} is the fact that the symplectic homology bases associated to two admissible decompositions are always related by an element of the group $\Gamma$. This fact is the content of Theorem~\ref{th1}, which is proven in Section~\ref{sect:Gam:SympBases}. The proof Theorem~\ref{th0} is then given in Section~\ref{sect:proofTh0}. In the Appendices, we give the proof of the fact that every surface in $\H(2)$ always admits parallelogram decompositions, and we give an explicit family of $\Gamma$-right cosets in $\Sp(4,\Z)$.

\bigskip

\textbf{ Acknowledgements :} The author is grateful to Universit\'e Paris Sud, Max Planck Institute for Mathematics, and Institute for Mathematical Sciences, NUS Singapore. He warmly thanks  François Labourie for the advices and stimulating discussions.

\section{Siegel upper half space and theta functions}\label{ThetSect}

In this section, we recall the definition and some basic properties of the theta functions. Our main references are \cite{RauFar}, and
\cite{FarKr}.

\subsection{ Siegel upper half space and Jacobian locus}

For any integer $g \geq 1$, the Siegel upper half space $\mathfrak{H}_g$ is the space of complex symmetric matrices $Z$ in $\mathfrak{M}(g,\C)$
such that $\mathrm{Im}(Z)$ is positive definite. This space is the quotient of the real symplectic group $\mathrm{Sp}(2g,\R)$ by the compact
subgroup $\mathrm{U}(g)$. The integral symplectic group $\mathrm{Sp}(2g,\Z)$ is a lattice in $\mathrm{Sp}(2g,\Z)$ which acts properly discontinuously on $\Hb_g$. The quotient $\mathcal{A}_g=\mathrm{Sp}(2g,\Z)\backslash\Hb_g$ is the moduli space of principally polarized Abelian varieties of
dimension $g$. In this paper, by the real symplectic group $\mathrm{Sp}(2g,\R)$ we mean the group of real $2g\times 2g$ matrices preserving the
symplectic form

$$\left(%
\begin{array}{ccc}
  J &  & 0 \\
    &\ddots & \\
  0 & & J \\
\end{array}%
\right), \text{ where } J=\left(%
\begin{array}{cc}
  0 & 1 \\
  -1 & 0\\
\end{array}%
\right).$$

Let $M$ be a Riemann surface of genus $g$, and $\{a_1,b_1,\dots, a_g,b_g\}$ be a symplectic basis of $H_1(M,\Z)$, that is

$$\langle a_i,a_j\rangle =\langle b_i,b_j\rangle =0, \text{ and }  \langle a_i,b_j\rangle=\delta_{ij},$$

\noindent where $\langle.,.\rangle$ is the intersection form of $H_1(M,\Z)$. There exist $g$ holomorphic $1$-forms $(\phi_1,\dots,\phi_g)$ on
$M$ uniquely determined by the following condition:

$$ \int_{a_j}\phi_i=\delta_{ij}.$$

\noindent The matrix $\Pi=(\pi_{ij})_{i,j=1,\dots,g}$, where $\pi_{ij}=\int_{b_j}\phi_i$, belongs to $\Hb_g$, and we then have a mapping from
the set of pairs $(M,\{a_1,b_1,\dots,a_g,b_g\})$ into $\Hb_g$. The image of this mapping is called the Jacobian locus denoted by
$\mathcal{J}_g$. 

For the case $g=2$, it is well known that the complement of $\mathcal{J}_2$ in $\Hb_2$ is a union of countably many copies of
$\Hb_1\times\Hb_1$, where $\Hb_1$ is the upper half plane $\Hb_1=\{z\in \C \; : \; \mathrm{Im} z >0 \}=\Hyp$ (see \cite{Mes}).

\subsection{Theta function}

Fix an integer $g\geq 1$, and let $\Hb_g$ be the { Siegel upper half space} of genus $g$. The Riemann's {\em theta function} is a complex value
function defined on $\C^g\times\Hb_g$ by the following formula

$$\theta(z,\sigma)=\sum_{N\in \Z^g} \exp(2\pi\imath(\frac{1}{2} \transpose{N}\sigma N+ ^tNz)).$$

\noindent The function $\theta$ is holomorphic on $\C^g\times\Hb_g$. We also consider functions defined on $\C^g\times \Hb_g$ by

$$\theta\left[%
\begin{array}{c}
  \epsilon \\
  \epsilon' \\
\end{array}%
\right](z,\sigma)=\sum_{N\in \Z^g} \exp\{ 2\pi\imath[\frac{1}{2}
 \transpose{(}N+\frac{\epsilon}{2})\sigma(N+\frac{\epsilon}{2})+ \transpose{(}N+\frac{\epsilon}{2})(z+\frac{\epsilon'}{2})]\}$$

\noindent where $\epsilon, \epsilon'$ are integer vectors. These functions are called first order theta functions with integer characteristic
$\left[\begin{array}{c}
  \epsilon \\
  \epsilon' \\
\end{array}\right]$.\\

\begin{proposition}\label{ThetProp}
The first order theta function with integer characteristic $\left[\begin{array}{c}
  \epsilon \\
  \epsilon' \\
\end{array}\right]$ has the following properties:

\begin{itemize}
\item[i)] $\theta\left[\begin{array}{c}
  \epsilon \\
  \epsilon' \\
\end{array}\right](z+e^{(k)},\sigma)=\exp 2\pi\imath \left[ \frac{\epsilon_k}{2}\right]\theta\left[\begin{array}{c}
  \epsilon \\
  \epsilon' \\
\end{array}\right](z,\sigma)$.\\

\item[ii)] $\theta\left[\begin{array}{c}
  \epsilon \\
  \epsilon' \\
\end{array}\right](z+\sigma^{(k)},\sigma)=\exp 2\pi \imath \left[ -z_k-\frac{\sigma_{kk}}{2}-\frac{\epsilon'_k}{2}\right] \theta\left[\begin{array}{c}
  \epsilon \\
  \epsilon' \\
\end{array}\right](z,\sigma)$.\\

\item[iii)] $\theta\left[\begin{array}{c}
  \epsilon \\
  \epsilon' \\
\end{array}\right](-z,\sigma)=\exp 2\pi\imath\left[\frac{\transpose{\epsilon}\epsilon'}{2}\right] \theta\left[\begin{array}{c}
  \epsilon \\
  \epsilon' \\
\end{array}\right](z,\sigma)$.\\

\item[iv)] $\theta\left[\begin{array}{c}
  \epsilon+2\nu \\
  \epsilon'+2\nu' \\
\end{array}\right](z,\sigma)=\exp 2\pi \imath\left[ \frac{\transpose{\epsilon}\nu'}{2}\right]\theta\left[\begin{array}{c}
  \epsilon \\
  \epsilon' \\
\end{array}\right](z,\sigma)$, with $\nu,\nu'\in \Z^g$,

\end{itemize}

\noindent where $e^{(k)}$ and $\sigma^{(k)}$ are the $k$-th column of the matrices $\Id_g$ and $\sigma$ respectively.

\end{proposition}

By Torelli Theorem, we know that a closed Riemann surface $M$ is uniquely determined by its Jacobian variety $J(M)$, or equivalently, by the
period matrix associated to a symplectic basis of $H_1(M,\Z)$. If $M$ is hyperelliptic, then we can get more information from the period matrix
by using theta functions. We have (cf. \cite{FarKr},  VII.4)

\begin{theorem} \label{ThetThm}
The branched points of the two sheeted representation of a hyperelliptic Riemann surface are holomorphic functions of the period matrix.
Furthermore, the hyperelliptic surface is completely determined by its period matrix.

\end{theorem}

\noindent To illustrate the ideas of the proof, we will indicate here below the method to compute some of the branched points, details of the
calculations can be found in \cite{FarKr}, VII.1, and VII.4.

Assume that $M$ is the two-sheeted branched cover of $\CP^1$ ramified above $\lbd_1,\dots,\lbd_{2g+2}$. Let $s_1,\dots,s_{g+1}$ be $g+1$ simple
arcs in $\CP^1$ such that the endpoints of $s_i$ are $\lbd_{2i-1},\lbd_{2i}$, for $k=i,\dots,g+1$, and $s_i\cap s_j=\vide$, for $i \neq j$.  We
can consider $M$ as the Riemann surface obtained by gluing two copies of $\CP^1$ slit along $s_1,\dots,s_{g+1}$. Let $z$ be the meromorphic
function on $M$ realizing the two-sheeted branched cover from $M$ to $\CP^1$ . Let $P_i$ denote $z^{-1}(\lbd_i), \; i=1,\dots,2g+2$, then
$\{P_1,\dots,P_{2g+2}\}$ is the set of Weierstrass points of $M$.

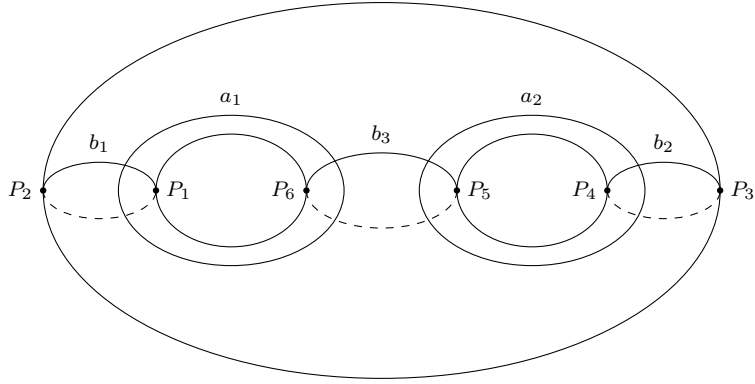
\begin{figure}[htb]
\begin{center}
\begin{tikzpicture}[scale=0.5]

\draw (0,0) ellipse (9 and 5); \draw (-4,0) ellipse (2 and 1.5) (4,0) ellipse (2 and 1.5);

\draw (-4,0) ellipse (3  and 2) (4,0) ellipse (3 and 2);

\draw (-6,0) arc (0:180:1.5 and 0.75) (9,0) arc (0:180:1.5 and 0.75); \draw[dashed] (-9,0) arc (180:360: 1.5 and 0.75) (6,0) arc (180:360:1.5
and 0.75);

\draw (2,0) arc (0:180:2 and 1); \draw[dashed] (-2,0) arc (180:360: 2 and 1);

\filldraw (-9,0) circle (2pt) (-6,0) circle (2pt) (-2,0) circle (2pt) (2,0) circle (2pt) (6,0) circle (2pt) (9,0) circle (2pt);

\tikzstyle{every node}=[font=\scriptsize]

\draw (-9,0) node[left] {$P_2$} (-6,0) node[right] {$P_1$} (-2,0) node[left] {$P_6$} (2,0) node[right] {$P_5$} (6,0) node[left] {$P_4$} (9,0)
node[right] {$P_3$};

\draw (-7.5,0.75) node[above] {$b_1$} (0,1) node[above] {$b_3$} (7.5,0.75) node[above] {$b_2$};

\draw (-4,2) node[above] {$a_1$} (4,2) node[above] {$a_2$};

\end{tikzpicture}
\end{center}
\caption{Symplectic basis on a hyperelliptic Riemann surface}
\label{fig:symp:basis:hyper}
\end{figure}

\noindent Using $\mathrm{PSL}(2,\C)$, we may assume that $\lbd_1=0,\lbd_2=1$, and $\lbd_{2g+2}=\infty$. The surface $M$ is then the curve defined by the equation

$$ w^2=z(z-1)\prod_{i=3}^{2g-1}(z-\lbd_i).$$

\noindent The function $z$ is then characterized, up to a non zero multiplicative constant, by the property that it has a double zero at $P_1$,
a double pole at $P_{2g+2}$, and it is holomorphic, and non zero elsewhere. First, we specify a symplectic  basis $\{a_1,b_1,\dots,a_g,b_g\}$ of $H_1(M,\Z)$ as follows:

\begin{itemize}

\item[$\bullet$] $b_k=z^{-1}(s_k), \; k=1,\dots,g+1$. By construction, $b_k$ is a simple closed curve containing $P_{2k-1}$,and $P_{2k}$, and
$b_k$ is reserved by the hyperelliptic involution.

\item[$\bullet$] Let $\alpha_k,\; k=1,\dots,g,$ be $g$ simple closed curves pairwise disjoint in $\CP^1$  satisfying

\begin{itemize}

\item[.] $\alpha_k$ intersects transversely $s_k$ and $s_{g+1}$,

\item[.] $\card \{\alpha_k\cap s_k\}=\card \{\alpha_k\cap s_{g+1}\}=1$,

\item[.] $\alpha_k\cap s_j=\vide$ if $j\notin \{k,g+1\}$,

\item[.] $\alpha_k\cap \{\lbd_1,\dots,\lbd_{2g+2}\}=\vide$.

\end{itemize}

Let $a_k$ be a connected component of $z^{-1}(\alpha_k)$. Note that $a_k$ and its image under the hyperelliptic  involution are disjoint.

\end{itemize}

\noindent It follows from the construction that the family $(a,b)=\{a_1,b_1,\dots,a_g,b_g\}$ is a symplectic basis of $H_1(M,\Z)$. Let
$\{\zeta_1,\dots,\zeta_g\}$ be the basis of $\mathscr{H}^1(M)$, the space of holomorphic $1$-form on $M$, dual to $\{a_1,\dots,a_g\}$, that is

$$\int_{a_j}\zeta_i=\delta_{ij}.$$

\noindent  Let $\Pi=(\pi_{ij})_{i,j=1,\dots,g}$ , with $\pi_{ij}=\int_{b_j}\zeta_i$, be the period matrix associated to $(a,b)$. By definition,
the Jacobian variety $J(M)$ of $M$ is the quotient $\C^g/\Lambda$, where $\Lambda$ is the lattice in $\C^g$ which is generated by the column
vectors of the $g\times 2g$ matrix $(\Id_g,\Pi)$. We denote by $e^{(k)}$, and $\pi^{(k)}, \; k=1,\dots,g,$  the $k$-th columns of $\Id_g$ and
$\Pi$ respectively. Let $\varphi: M \lra J(M)$ be the map defined by

$$\varphi(P)=(\int_{P_1}^P\zeta_1,\dots,\int_{P_1}^P\zeta_g)\in \C/\Lambda,$$

\noindent where the integrals are taken along any path joining $P_1$ to $P$. From the construction of the basis $(a,b)$, one can compute
explicitly  $\varphi(P_i), \; i=1,\dots,2g+2$, as functions of $(e^{(1)},\dots,e^{(g},\pi^{(1)},\dots,\pi^{(g)})$. Namely, we have (see Figure~\ref{fig:symp:basis:hyper})

\begin{itemize}
\item[$\bullet$] $\varphi(P_1)=0$.

\item[$\bullet$] $\varphi(P_2)=\frac{1}{2}\pi^{(1)}$.

\item[$\bullet$] $\varphi(P_3)= \frac{1}{2}(\pi^{(1)}+e^{(1)}+e^{(2)})$.

\item[$\bullet$] $\dots$

\item[$\bullet$] $\varphi(P_{2k+1})=\frac{1}{2}(\pi^{(1)}+\dots+\pi^{(k)}+e^{(1)}+\dots+e^{(k+1)})$.

\item[$\bullet$] $\varphi(P_{2k+2})=\frac{1}{2}(\pi^{(1)}+\dots+\pi^{(k+1)}+e^{(1)}+\dots+e^{(k+1)})$.

\item[$\bullet$] $\dots$

\item[$\bullet$] $\varphi(P_{2g+1})=\frac{1}{2}(\pi^{(1)}+\dots+\pi^{(g)}+e^{(1)})$.

\item[$\bullet$] $\varphi(P_{2g+2})=\frac{1}{2}e^{(1)}$.\\

\end{itemize}

Since $\Pi \in \Hb_g$, we can now consider the first order theta functions with integer characteristic $\theta\left[\begin{array}{c} \epsilon \\
\epsilon' \\ \end{array}\right](z,\Pi)$. From Proposition \ref{ThetProp}, we see that for any $\nu \in \Lambda$, $\theta\left[\begin{array}{c}
\epsilon \\ \epsilon' \\ \end{array}\right](z+\nu,\Pi)$ differs from $\theta\left[\begin{array}{c} \epsilon \\ \epsilon' \\
\end{array}\right](z,\Pi)$ by a multiplicative factor. It turns out that the multiplicative behavior of the theta functions is so that

$$f(P):=\frac{\theta^2\left[\begin{array}{ccccc} 1 & 0 & 0 & \dots &0\\ 1 & 1 & 0 & \dots & 0\\ \end{array}\right](\varphi(P),\Pi)}
{\theta^2\left[\begin{array}{ccccc} 1 & 0 &  0 &\dots& 0 \\ 0 & 1 & 0 & \dots & 0\\ \end{array}\right](\varphi(P),\Pi)}$$

\noindent is a meromorphic function on $M$ with divisor $P_1^2P_{2g+2}^{-2}$. Hence

$$f=cz, \text{ where } c\in \C^*.$$

\noindent The constant $c$ can be valuated at $P_2$ since we have  $f(P_2)=cz(P_2)=c$. Using the fact that $\varphi(P_2)=\frac{1}{2}\pi^{(1)}$,
we have

$$c=f(P_2)=\frac{\theta^2\left[\begin{array}{ccccc} 1 & 0 & 0 & \dots & 0 \\ 1 & 1 & 0 &\dots& 0 \\ \end{array} \right](\frac{1}{2}\pi^{(1)},\Pi)}
{\theta^2\left[\begin{array}{ccccc} 1 & 0 & 0 & \dots & 0 \\ 0 & 1 & 0 &\dots& 0 \\ \end{array} \right](\frac{1}{2}\pi^{(1)},\Pi)}.$$

\noindent It follows that

\begin{equation}\label{ThetEq}
z(P)= \frac{ \theta^2\left[\begin{array}{ccccc} 1 & 0 & 0 & \dots & 0 \\ 0 & 1 & 0 &\dots& 0 \\ \end{array} \right](\frac{1}{2}\pi^{(1)},\Pi) \theta^2 \left[\begin{array}{ccccc} 1 & 0 & 0 & \dots & 0 \\ 1 & 1 & 0 &\dots& 0 \\ \end{array} \right](\varphi(P),\Pi)}
{\theta^2\left[\begin{array}{ccccc} 1 & 0 & 0 & \dots & 0 \\ 1 & 1 & 0 &\dots& 0 \\ \end{array} \right](\frac{1}{2}\pi^{(1)},\Pi) \theta^2 \left[\begin{array}{ccccc} 1 & 0 & 0 & \dots & 0 \\ 0 & 1 & 0 &\dots& 0 \\ \end{array} \right](\varphi(P),\Pi)}
\end{equation}

\noindent By setting $P=P_j$, for $j=3,\dots,2g+1$, we get a formula for $\lbd_j$. Note that this formula is only useful for $j=3,4,6,..,2g$,
for other values of $j$, both the numerator and denominator vanish. By replacing $P_1$ by another Weierstrass point in the definition of
$\varphi$, we can get similar formulae for the $\lbd_j$ which can not be computed directly from (\ref{ThetEq}).

\section{The group $\Gamma$}\label{GrpGamSect}

In this section, we define the group $\Gamma$, and prove some of its properties.

\subsection{Construction of surfaces in $\H(2)$ by gluing parallelograms}

Any parallelogram in $\R^2$ is determined by a pair of vectors in $\R^2$. In this section, we will represent any parallelogram by a pair of
complex numbers $(z_1,z_2)$ such that $\DS{\mathrm{Im}(\bar{z}_1z_2)=\frac{-\imath}{2}(\bar{z}_1z_2-z_1\bar{z}_2)>0}$. Note that by this
convention, the pairs $(z_1,z_2), (z_2,-z_1), (-z_1,-z_2), (-z_2,z_1)$ represent the same parallelogram.

Let $\mathscr{P}^+$ denote the set $\DS{\{(z_1,z_2)\in \C^2\; : \frac{-\imath}{2}(\bar{z}_1z_2-z_1\bar{z}_2)>0\}}$. Given $4$ complex numbers
$(z_1,\dots,z_4)$ such that $(z_1,z_2), (z_2,z_3), (z_3,z_4)$ all belong to $\mathscr{P}^+$, let $A,B,C$ denote the parallelograms determined by
the pairs $(z_1,z_2), (z_2,z_3),(z_3,z_4)$ respectively. We can construct a translation surface in $\H(2)$ from $A,B,C$ as follows

\begin{itemize}

\item[$\bullet$] Glue two sides of $A$ corresponding to $z_1$ together,

\item[$\bullet$] Glue two sides of $A$ corresponding to $z_2$ to two sides of $B$ also corresponding to $z_2$,

\item[$\bullet$] Glue two sides of $B$ corresponding to $z_3$ to two sides of $C$ also corresponding to $z_3$,

\item[$\bullet$] Glue two sides of $C$ corresponding to $z_4$ together.

\end{itemize}

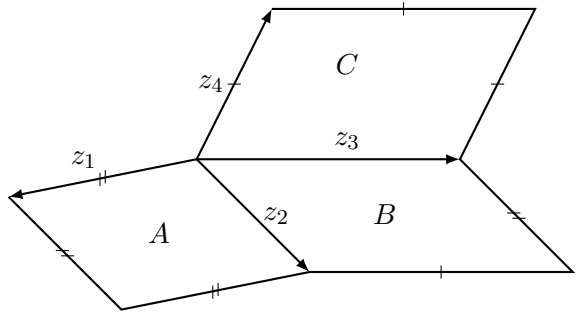
\begin{figure}[htb]
\begin{center}
\begin{tikzpicture}[>=latex, scale=0.5]
\draw[->, thick] (0,0)-- (2,4); \draw[->, thick] (0,0) -- (7,0); \draw[->, thick] (0,0) -- (3,-3); \draw[->, thick] (0,0) -- (-5,-1);
\draw[thick] (2,4) -- (9,4) -- (7,0) -- (10,-3) -- (3,-3) -- (-2,-4) -- (-5,-1);

\draw (1,2) node[left] {$z_4$}; \draw (4,0) node[above] {$z_3$}; \draw (1.5,-1.5) node[right] {$z_2$}; \draw (-3,-0.5) node[above] {$z_1$};

\draw (4,2) node[above] {$C$}; \draw (5,-2) node[above] {$B$}; \draw (-1,-2.5) node[above] {$A$};


\draw (1,2)  +(-5pt,0)-- +(5pt,0); \draw (8,2) +(-5pt,0) -- +(5pt,0); \draw (-3.5,-2.5) +(-7pt,2pt) -- +(3pt,2pt) +(-3pt,-2pt) -- +(7pt,-2pt);
\draw (8.5,-1.5) +(-7pt,2pt) -- +(3pt,2pt) +(-3pt,-2pt) -- +(7pt,-2pt); \draw (5.5,4) +(0,-5pt) -- +(0,5pt); \draw (6.5,-3) +(0,-5pt) --
+(0,5pt); \draw (-2.5,-0.5) +(-2pt, -6pt) -- +(-2pt,4pt) +(2pt,-4pt) -- +(2pt,6pt); \draw (0.5,-3.5) +(-2pt, -6pt) -- +(-2pt,4pt) +(2pt,-4pt) --
+(2pt,6pt);

\end{tikzpicture}
\end{center}
\caption{Construction of surfaces in $\H(2)$ from three parallelograms.}
\label{fig:3par:construction}
\end{figure}

It is easy to check that the surface $M$ obtained from  this construction is of genus $2$, equipped with a flat metric structure with a single cone
singularity, which arises from the identification of all the vertices of $A,B,C$, we denote this point by $W$. Since all the gluings are realized by translations, $M$ is a translation surface. We also get naturally a holomorphic $1$-form on $M$, considered as a Riemann surface, defined as follows: since translations of $\R^2$ preserve the holomorphic $1$-form $dz$ is preserved, the restrictions of $dz$ into the parallelograms $A,B,C$ are compatible with the gluings, and give rise to a holomorphic $1$-form $\omega$ on $M$ with only one zero at $W$, which is necessarily of order two. Clearly, $(M,\omega)\in \H(2)$. 

We know that $M$ is a hyperelliptic surface, it is easy to visualize the hyperelliptic involution of $M$ from its construction by gluing
$A,B,C$. For each of the parallelograms $A,B,C$, consider the refection through its center, one can easily  check that these reflections agree
with the gluing on the boundary of $A,B,C$. Thus, we have a conformal automorphism $\tau$ of $M$. One can check that $\tau^2=\Id$, and the
action of $\tau$ on $H_1(M,\Z)$ is given by $-\Id$, therefore $\tau$ must be the hyperelliptic involution of $M$. We can also determine without
difficulty the $6$ fixed points of $\tau$, which are the Weierstrass points of $M$: two of which are contained in $A$, two in $C$, one in the
interior of $B$, and the last one is $W$.

\subsection{Parallelogram decompositions and the group $\Gamma$}\label{ParDecSect}

Recall that, on translation surfaces $(M,\omega)$, a {\em saddle connection} is a geodesic segment (with respect to the flat metric) which joins singularity to singularity, the endpoints of a saddle connection may coincide. A {\em cylinder} $C$ is a subset of $M$ which is isometric to $\R\times]0,h[/\Z$, where the action of $\Z$ is generated by $(x,y) \mapsto (x+\ell,y)$, with $\ell>0$, and maximal with respect to this property (that is, $C$ can not be embedded into a larger subset $C'$ isometric to $\R\times ]0,h'[/ \Z$, with $h'>h$). In other words, $C$ is the union of all simple closed geodesics in some free homotopy class. The construction of translation surfaces in $\H(2)$ by gluing parallelograms with the model presented above suggests the following

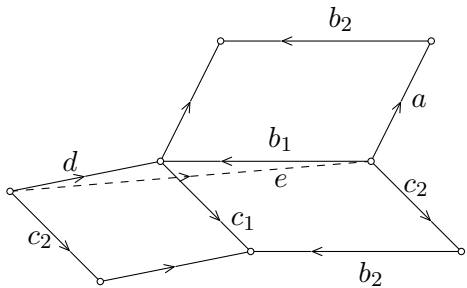
\begin{figure}[htb]
\begin{center}
\begin{tikzpicture}[scale=0.4]
\draw[->, >=angle 45] (0,0)-- (1,2); \draw (1,2) -- (2,4);
\draw[->, >= angle 45] (7,0) -- (8,2); \draw (8,2) -- (9,4);

\draw[->, >=angle 45] (9,4) -- (4,4); \draw (4,4) -- (2,4);
\draw[->, >= angle 45] (7,0) -- (2,0); \draw (2,0) -- (0,0);
\draw[->, >=angle 45] (10,-3) -- (5,-3); \draw (5,-3) -- (3,-3);

\draw[->, >=angle 45] (-5,-1) -- (-3,-3); \draw (-3,-3) -- (-2,-4);
\draw[->, >= angle 45] (0,0) -- (2,-2); \draw (2,-2) -- (3,-3);
\draw[->, >=angle 45] (7,0) -- (9,-2); \draw (9,-2) -- (10,-3);

\draw[->,>=angle  45] (-5,-1) -- (-2.5,-0.5); \draw (-2.5,-0.5) -- (0,0);
\draw[->, >=angle 45] (-2,-4) -- (0.5,-3.5); \draw (0.5,-3.5) -- (3,-3);
\draw[dashed,->,>=angle 45] (-5,-1) -- (1,-0.5); \draw[dashed] (1,-0.5) -- (7,0); 

\foreach \x in {(-5,-1), (-2,-4), (0,0), (2,4), (3,-3), (7,0), (9,4), (10,-3)} \filldraw[fill=white] \x circle (3pt);

\draw (8,2) node[right] {$a$} (6,4) node[above] {$b_2$} (4,0) node[above] {$b_1$} (7,-3) node[below] {$b_2$} (8.5,-1.5) node[above] {$c_2$} (2,-2) node[right] {$c_1$} (-4,-2) node[below] {$c_2$} (-3,0) node {$d$} (4,-1.2) node[above] {$e$};

\end{tikzpicture}
\end{center}
\caption{Parallelogram decomposition of surfaces in $\H(2)$.}
\label{fig:ParDecConfig}

\end{figure}

\begin{definition}\label{PDcpDef}
Let $(M,\omega)$ be a pair in $\H(2)$. A {\em parallelogram decomposition } of $(M,\omega)$ is a family of six oriented saddle
connections $\{a,b_1,b_2,c_1,c_2,d\}$ verifying the following conditions

\begin{itemize}

\item[$\bullet$] The intersection of any pair of saddle connections in this family is the only zero of $\omega$,

\item[$\bullet$] $b_1\cup b_2$ (resp. $c_1\cup c_2$) is the boundary of a cylinder which contains $a$ (resp. $d$),

\item[$\bullet$] The complement of $a\cup b_1\cup b_2 \cup c_1 \cup c_2 \cup d$ has three components, each of which is isometric to an open parallelograms in $\R^2$,

\item[$\bullet$] The orientations of the saddle connections in this family are chosen so that  $\langle a, b_1 \rangle = \langle b_1, c_1 \rangle =\langle c_1, d \rangle =1$,  $a$ goes from $b_1$ to $b_2$, and $d$ goes from $c_2$ to $c_1$ (see Figure~\ref{fig:ParDecConfig}).

\end{itemize}

\end{definition}

If $\{a,b_1,b_2,c_1,c_2,d\}$ is a parallelogram decomposition of $(M,\omega)$ then $(a,b_1,c_1,d)$ is a non-symplectic basis of $H_1(M,\Z)$.
Let $b$ (resp. $c$) be a simple closed curve in the cylinder bounded by $b_1$ and $b_2$ (resp. by $c_1$ and $c_2$). Let $e$ be a simple closed
curved in the free homotopy class of the  closed curve $d*(-b_1)$, that is $e=d-b$ in $H_1(M,\Z)$, then $(a,b,c,e)$ is a symplectic basis of $H_1(M,\Z)$, we will call it the {\em symplectic basis associated to } $\mathcal{D}$.

Given a surface  $(M,\omega)$ in $\H(2)$ which is obtained from three parallelograms $A=(z_1,z_2),B=(z_2,z_3),C=(z_3,z_4)$ as in the previous
section, let $\tilde{A},\tilde{B},\tilde{C}$ be the subsets of $M$ which correspond to $A,B,C$ respectively. By construction, $\tilde{A}$ and
$\tilde{C}$ are two cylinders, while $\tilde{B}$ is an embedded parallelogram. Let $a$ (resp. $d$) denote the saddle connection in $\tilde{A}$ (resp.
$\tilde{C})$ which corresponds to the $z_1$ sides of $A$ (resp. $z_4$ sides of $C$). Let $b_1,b_2$ (resp. $c_1,c_2$) denote the boundary
components of $\tilde{A}$ (resp. $\tilde{C}$). We choose the orientations of $a,b_1,b_2,c_1,c_2,d$ so that

$$\int_{a}\omega= z_1, \; \int_{b_1}\omega=\int_{b_2}\omega=z_2, \; \int_{c_1}\omega=\int_{c_2}\omega=z_3, \; \int_d\omega=z_4.$$

\noindent We also choose the numbering of $(b_1,b_2)$ (resp. $(c_1,c_2)$) so that the orientation of $a$ goes from $b_1$ to $b_2$, and the
orientation of $d$ goes from $c_2$ to $c_1$. By definition $\{a,b_1,b_2,c_1,c_2,d\}$ is a parallelogram decomposition of $(M,\omega)$.

A surface $(M,\omega)$ in $\H(2)$ always admits a parallelogram decompositions see Appendix Section~\ref{app:sect:ParDecomp:Existence}. The following operations allow us to get other decompositions from $\mathcal{D}=\{a,b_1,b_2,c_1,c_2,d\}$.

\begin{itemize}

\item[1.] \underline{\bf The $T$ move:} let $a'$ be the saddle connection  which is obtained from $a$ by a Dehn twist in $\tilde{A}$, then
$\mathcal{D}'=\{a',b_1,b_2,c_1,c_2,d\}$ is another parallelogram decomposition of $(M,\omega)$ (see Figure~\ref{fig:Tmove}).

\begin{figure}[htb]
\begin{center}
\begin{tikzpicture}[scale=0.3]


\draw[>=angle 45, ->] (-4,3) -- (-2,3.5); \draw[->, >=angle 45] (-4,-1) -- (-2,-0.5); \draw[->, >=angle 45] (0,4) -- (0,1); \draw[->, >=angle
45] (-4,3) -- (-4,0); \draw[->, >=angle 45] (6,4) -- (2,4); \draw[->, >=angle 45] (4,7) -- (0,7); \draw[->, >=angle 45] (6,0) -- (2,0);
\draw[->, >=angle 45] (6,4) -- (5,5.5); \draw[->, >=angle 45] (0,4) -- (-1,5.5); \draw[->, >=angle 45] (6,4) -- (6,1);

\draw (-2,3.5) -- (0,4) -- (2,4) (-2,-0.5) -- (0,0) -- (2,0) (-4,-1) -- (-4,0) (6,1) -- (6,0) (0,0) -- (0,1) (5,5.5) -- (4,7) (-2,7) -- (-1,5.5)
(-2,7) -- (0,7);

\draw (5,5.5) node[right] {$a$} (3,4) node[above] {$b_1$} (1,7) node[above] {$b_2$} (3,0) node[below] {$b_2$} (0,2) node[right] {$c_1$} (-4,1)
node[left] {$c_2$} (6,2) node[left] {$c_2$} (-2,3.5) node[above] {$d$};

\draw[->, >=angle 45] (11,3) -- (13,3.5); \draw[->, >=angle 45] (11,3) -- (11,0); \draw[->, >=angle 45] (11,-1) -- (13,-0.5); \draw[->, >=angle
45] (15,4) -- (17,5.5); \draw[->, >=angle 45] (15,4) -- (15,1); \draw[->, >=angle 45] (21,0) -- (17,0); \draw[->, >=angle 45] (21,4) --
(23,5.5); \draw[->, >=angle 45] (21,4) -- (17,4); \draw[->, >=angle 45] (21,4) -- (21,1); \draw[->, >=angle 45] (25,7) -- (21,7);

\draw (11,0) -- (11,-1) (13,-0.5) -- (15,0) -- (17,0) (15,0) -- (15,1) (13,3.5) -- (15,4) -- (17,4) (21,1) -- (21,0) (17,5.5) -- (19,7) --
(21,7) (23,5.5) -- (25,7);

\draw (23,5.5) node[right] {$a'$} (18,4) node[above] {$b_1$} (22,7) node[above] {$b_2$} (18,0) node[below] {$b_2$} (15,2) node[right] {$c_1$}
(11,1) node[left] {$c_2$} (21,2) node[right] {$c_2$} (13,3.5) node[above] {$d$};

\draw[dashed] (0,4) -- (4,7);

\draw[->, >= angle 90, thick, red, decorate, decoration={snake, amplitude=3, segment length=6, post length=3}] (7,4) -- (10,4);

\end{tikzpicture}
\end{center}
\caption{The $T$ move.}
\label{fig:Tmove}
\end{figure}
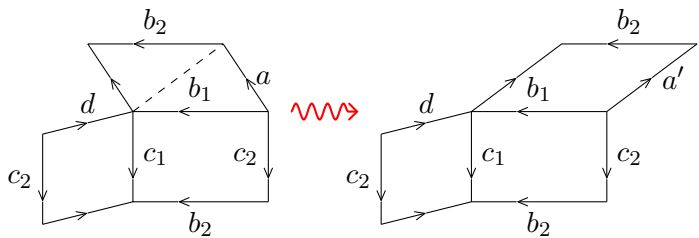

\item[2.] \underline{\bf The $S$ move:} If $\{a,b_1,b_2,c_1,c_2,d\}$ is a parallelogram decomposition of $(M,\omega)$ then\\ $\De'=
\{d,-c_2,-c_1, b_1,b_2,-a\}$ is also a parallelogram decomposition of $(M,\omega)$, the minus sign designate the same saddle connection with the
inverse orientation (see Figure~\ref{fig:Smove}).

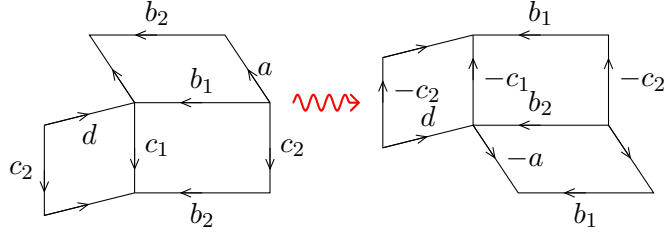
\begin{figure}
\begin{center}
\begin{tikzpicture}[scale=0.3, >=angle 45]


\draw (-4,3) -- (-4,-1) -- (0,0) -- (6,0) -- (6,4) -- (4,7) -- (-2,7) -- (0,4) -- (6,4); \draw (-4,3) -- (0,4) -- (0,0);

\draw[->] (-4,3) -- (-2,3.5); \draw[->] (-4,-1) -- (-2,-0.5); \draw[->] (-4,1) -- (-4,0); \draw[->] (0,2) -- (0,1); \draw[->] (6,2) -- (6,1);
\draw[->] (3,0) -- (2,0); \draw[->] (3,4) -- (2,4); \draw[->] (1,7) -- (0,7); \draw[->] (6,4) -- (5,5.5); \draw[->] (0,4) -- (-1,5.5);


\draw (11,2) -- (11,6) -- (15,7) -- (21,7) -- (21,3) -- (23,0) -- (17,0) -- (15,3) -- cycle; \draw (15,7) -- (15,3) -- (21,3);

\draw[->] (11,4) -- (11,5); \draw[->] (15,5) -- (15,6); \draw[->] (21,5) -- (21,6); \draw[->] (11,2) -- (13,2.5); \draw[->] (11,6) --(13,6.5);
\draw[->] (18,7) -- (17,7); \draw[->] (18,3) -- (17,3); \draw[->] (20,0) -- (19,0); \draw[->] (15,3) -- (16,1.5); \draw[->] (21,3) -- (22,1.5);

\draw[->, >=angle 90, thick, red, decorate, decoration={snake, amplitude=3, segment length=6, post length=3}] (7,4) -- (10,4);

\draw (5,5.5) node[right] {$a$}; \draw (3,4) node[above] {$b_1$}; \draw (0,2) node[right] {$c_1$}; \draw (-2,3.5) node[below] {$d$}; \draw (1,7)
node[above] {$b_2$}; \draw (3,0) node[below] {$b_2$}; \draw (-4,1) node[left] {$c_2$} (6,2) node[right] {$c_2$};

\draw (13,2.5) node[above] {$d$}; \draw (15,5) node[right] {$-c_1$}; \draw (18,3) node[above] {$b_2$}; \draw (16,1.5) node[right] {$-a$}; \draw
(18,7) node[above] {$b_1$}; \draw (20,0) node[below] {$b_1$}; \draw (11,4.5) node[right] {$-c_2$} (21,5) node[right] {$-c_2$};

\end{tikzpicture}

\end{center}
\caption{The $S$ move.}
\label{fig:Smove}
\end{figure}

\item[3.] \underline{\bf The $R$ move} Cut $M$ along the saddle connections $b_1,b_2$ and $d$, we then obtain two cylinders, one of which is
$\tilde{A}$, we denote the other by $\tilde{A}'$. Let $c'_1,c'_2$ denote the images of $c_1,c_2$ respectively under a Dehn twist in
$\tilde{A}'$. Assume that $c'_1,c'_2$ can be made  into saddle connections in  $\tilde{A}'$, then $\{a,b_1,b_2,c'_1,c'_2,d\}$ is a parallelogram
decomposition of $(M,\omega)$ (see Figure~\ref{fig:Rmove}).

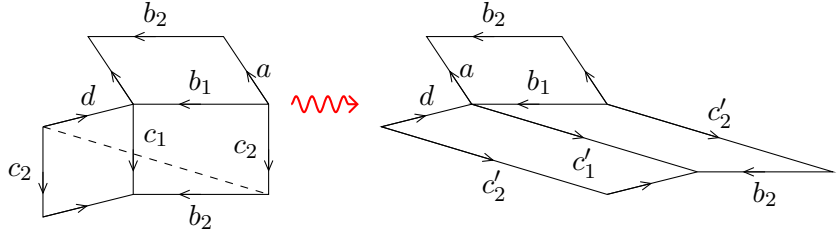
\begin{figure}[htb]
\begin{center}
\begin{tikzpicture}[scale=0.3, >=angle 45]


\draw (-4,3) -- (-4,-1) -- (0,0) -- (6,0) -- (6,4) -- (4,7) -- (-2,7) -- (0,4) -- (6,4); \draw (-4,3) -- (0,4) -- (0,0);

\draw[->] (-4,3) -- (-2,3.5); \draw[->] (-4,-1) -- (-2,-0.5); \draw[->] (-4,1) -- (-4,0); \draw[->] (0,2) -- (0,1); \draw[->] (6,2) -- (6,1);
\draw[->] (3,0) -- (2,0); \draw[->] (3,4) -- (2,4); \draw[->] (1,7) -- (0,7); \draw[->] (6,4) -- (5,5.5); \draw[->] (0,4) -- (-1,5.5);

\draw[dashed] (-4,3) -- (6,0);

\draw (5,5.5) node[right] {$a$} (3,4) node[above] {$b_1$} (1,7) node[above] {$b_2$} (3,0) node[below] {$b_2$} (0,2.5) node[right] {$c_1$} (-4,1)
node[left] {$c_2$} (6,2) node[left] {$c_2$} (-2,3.5) node[above] {$d$};


\draw (11,3) -- (21,0) -- (25,1) -- (31,1) -- (21,4) -- (19,7) -- (13,7) -- (15,4) -- cycle; \draw (25,1) -- (15,4) -- (21,4);

\draw[->] (15,4) -- (14,5.5); \draw[->] (21,4) -- (20,5.5); \draw[->] (16,7) -- (15,7); \draw[->] (18,4) -- (17,4); \draw[->] (28,1) --(27,1);
\draw[->] (11,3) -- (13,3.5); \draw[->] (21,0) -- (23,0.5); \draw[->] (11,3) -- (16,1.5); \draw[->] (15,4) -- (20,2.5); \draw[->] (21,4) --
(26,2.5);

\draw (14,5.5) node[right] {$a$} (18,4) node[above] {$b_1$} (16,7) node[above] {$b_2$} (28,1) node[below] {$b_2$} (20,2.5) node[below] {$c'_1$}
(16,1.5) node[below] {$c'_2$} (26,2.5) node[above] {$c'_2$} (13,3.5) node[above] {$d$};

\draw[->, >=angle 90, thick, red, decorate, decoration={snake, amplitude=3, segment length=6, post length=3}] (7,4) -- (10,4);

\end{tikzpicture}
\end{center}
\caption{The $R$ move.}
\label{fig:Rmove}
\end{figure}

\end{itemize}

Let $\mathcal{D}'$ be another parallelogram decomposition of $(M,\omega)$ which is obtained from $\De$ by one the moves presented above. Let $(a,b,c,e)$ (resp. $(a',b',c',e'))$ be the symplectic basis of $H_1(M,\Z)$ associated to $\mathcal{D}$ (resp. $\mathcal{D}')$. The following lemma follows directly from the definitions.

\begin{lemma}\label{lm:elmntmoves:matrices}
\begin{itemize}

\item[a)] If $\mathcal{D}'$ is obtained from $\mathcal{D}$ by a $T$ move then

$$\left(%
\begin{array}{c}
  a' \\
  b'\\
  c' \\
  e' \\
\end{array}%
\right)=\left(%
\begin{array}{cccc}
  1 & \pm 1 & 0 & 0 \\
  0 & 1 & 0 & 0 \\
  0 & 0 & 1 & 0 \\
  0 & 0 & 0 & 1 \\
\end{array}%
\right) \cdot \left(%
\begin{array}{c}
  a \\
  b \\
  c \\
  e \\
\end{array}%
\right).$$

\item[b)] If $\mathcal{D}'$ is obtained from $\mathcal{D}$ by an $S$ move then

$$\left(%
\begin{array}{c}
  a' \\
  b'\\
  c' \\
  e' \\
\end{array}%
\right)= \left(%
\begin{array}{cccc}
  0 & 1 & 0 & 1 \\
  0 & 0 & -1 & 0 \\
  0 & 1 & 0 & 0 \\
  -1 & 0 & 1 & 0 \\
\end{array}%
\right) \cdot \left(%
\begin{array}{c}
  a \\
  b \\
  c \\
  e \\
\end{array}%
\right).$$

\item[c)] If $\mathcal{D}'$ is obtained from $\mathcal{D}$ by an $R$ move then

$$\left(%
\begin{array}{c}
  a' \\
  b'\\
  c' \\
  e' \\
\end{array}%
\right)= \left(%
\begin{array}{cccc}
  1 & 0 & 0 & 0 \\
  0 & 1 & 0 & 0 \\
  0 & 0 & 1 & \pm 1 \\
  0 & 0 & 0 & 1 \\
\end{array}%
\right) \cdot \left(%
\begin{array}{c}
  a \\
  b \\
  c \\
  e \\
\end{array}%
\right).$$

\end{itemize}
\end{lemma}

We denote by $T,S,R$ the matrices of basis change corresponding to the moves $T,S,R$ respectively, and by $\Gamma$ the subgroup of $\mathrm{Sp}(4,\Z)$ generated by those matrices.


\subsection{ Properties of $\Gamma$}

\begin{lemma}\label{GamPropLm}
We have
\begin{itemize}

\item[(i)] $S^2=-\Id_4$.

\item[(ii)] $\left(%
\begin{array}{cc}
  \Id_2 & 0 \\
  0 & \mathrm{SL}(2,\Z)
\end{array}%
\right) \subset \Gamma$.

\item[(iii)] $\Gamma \varsubsetneq \mathrm{Sp}(4,\Z)$.

\item[(iv)] $\Gamma$ is not a normal subgroup of $\mathrm{Sp}(4,\Z)$.

\end{itemize}
\end{lemma}

\dem

\begin{itemize}

\item[(i)] follows from direct calculation.

\item[(ii)] We have

$$ S^{-1}TS=\left(%
\begin{array}{cccc}
  1 & 0 & 0 & 0 \\
  0 & 1 & 0 & 0 \\
  0 & 0 & 1 & 0 \\
  0 & 0 & -1 & 1 \\
\end{array}%
\right).$$

Since $\mathrm{SL}(2,\Z)$ is generated by $\DS{ \left(%
\begin{array}{cc}
  1 & 1 \\
  0 & 1 \\
\end{array}%
\right), \left(%
\begin{array}{cc}
  1 & 0 \\
  -1 & 1 \\
\end{array}%
\right)} $, we deduce that $\DS{\left(\begin{array}{cc} \Id_2 & 0 \\ 0 & \SL(2,\Z)\\ \end{array} \right)}$ is contained in $\Gamma$.

\item[(iii)] The group $\mathrm{Sp}(4,\Z)$ acts transitively on $(\Z/2\Z)^4\setminus \{0\}$, but $\Gamma$ has two orbits: $\mathscr{O}_1$ containing
$\mathbf{e}_1=(1,0,0,0)$, and $\mathscr{O}_2$ containing $\mathbf{e}_2=(0,1,0,0)$. As a matter of fact, we have

$$\mathscr{O}_1=\left\{\begin{array}{ccc} (1,0,0,0) & (1,1,0,0) &  \\ (0,1,1,0) & (0,1,0,1) & (0,1,1,1) \end{array} \right\}$$

\noindent and

$$\mathscr{O}_2=\left\{ \begin{array}{ccc} (0,1,0,0) &  & \\ (1,0,1,0) & (1,0,0,1) & (1,0,1,1) \\ (1,1,1,0) & (1,1,0,1) & (1,1,1,1) \\ (0,0,1,0) & (0,0,0,1) & (0,0,1,1) \end{array}\right\}$$

\noindent Here, we consider the action of $\Sp(4,\Z)$ and $\Gamma$ on $(\Z/2\Z)^4$ by right multiplication. 


\item[(iv)] Let $T'=\left(%
\begin{array}{cccc}
  1 & 0 & 0 & 0 \\
  1 & 1 & 0 & 0 \\
  0 & 0 & 1 & 0 \\
  0 & 0 & 0 & 1 \\
\end{array}%
\right)$. Remark that $T'$ does not belong to $\Gamma$ since it sends $\mathbf{e}_2$ to an element in $\mathscr{O}_1$. We have

$$ {T'}^{-1}TT'=\left(%
\begin{array}{cccc}
  2 & 1 & 0 & 0 \\
  -1 & 0 & 0 & 0 \\
  0 & 0 & 1 & 0 \\
  0 & 0 & 0 & 1 \\
\end{array}%
\right).$$

Observe that ${T'}^{-1}T T'$ sends $\mathbf{e}_1$ to $\mathbf{e}_2$, thus it does not belong to $\Gamma$. We can then conclude that
$\Gamma$ is not a normal subgroup of $\mathrm{Sp}(4,\Z)$.

\end{itemize}

\carre

\begin{lemma}\label{SpGenLm}
Set

$$U=\left(%
\begin{array}{cccc}
  0 & 0 & 0 & -1 \\
  0 & 0 & 1 & 0 \\
  0 & -1 & 0 & 0 \\
  1 & 0 & 0 & 0\\
\end{array}%
\right),$$

\noindent then the integral symplectic group $\mathrm{Sp}(4,\Z)$ is generated by $T,S,R$, and $U$.

\end{lemma}

\dem see Appendix~\ref{prfSpGenLm}.\carre

\section{Admissible decomposition}\label{AdmDecSect}

\subsection{Definition}

We have seen that for any surface $(M,\omega)$ in $\H(2)$, given a parallelogram decomposition, other decompositions of $(M,\omega)$ can be obtained by the elementary moves $T,S$, and $R$. However, while $T$ and $S$ are always realizable, the $R$ move is not,  it is only realizable when $d*(-b_1)*c_2$ is homotopic to a saddle connection. In this section, we  enlarge the set of decompositions  so that the three elementary moves are always realizable. This leads us to the notion of {\em admissible decomposition} with respect to the pair $(M,W)$. The definition of admissible decomposition is inspired from parallelogram decomposition, and based on the action of the hyperelliptic involution of $M$.

Throughout this section, $M$ is a fixed closed Riemann surface of genus  two, $W$ is a Weierstrass point of $M$. Let $\tau$ denote the
hyperelliptic involution of $M$. For any closed curve $\gamma$ with base point $W$, we denote by $[\gamma]$ the homotopy class of $\gamma$ in
$\pi_1(M,W)$.

\begin{definition}[Admissible decomposition]\label{AdmDecDef}

Let $\{a,b_1,b_2,c_1,c_2,d\}$ be six oriented simple closed curve with base point at $W$. We say that $\{a,b_1,b_2,c_1,c_2,d\}$ is an {\em
admissible decomposition} for the pair $(M,W)$ if

\begin{itemize}

\item[$\bullet$] The intersection of any pair of curves in this family is $\{W\}$.

\item[$\bullet$] $\tau(a)=-a, \tau(d)=-d$.

\item[$\bullet$] $\tau(b_1)=-b_2, \; \tau(c_1)=-c_2$.

\item[$\bullet$] $a\setminus\{W\}$ is contained in an open annulus $A_b$ bounded by $b_1$ and $b_2$.

\item[$\bullet$] $d\setminus \{W\}$ is contained in an open annulus $A_c$ bounded by $c_1$ and $c_2$.

\item[$\bullet$] $M\setminus \overline{A}_b \cup \overline{A}_c$ is homeomorphic to an open disk.

\end{itemize}


The orientations of $\{a,b_1,b_2,c_1,c_2,d\}$ are chosen so that

\begin{itemize}

\item[$\bullet$] $\langle a,b_1\rangle =\langle a,b_2\rangle =\langle c_1,d \rangle= \langle c_2,d \rangle =1$, where $\langle .,. \rangle$ is
the intersection form of $H_1(M,\Z)$.

\item[$\bullet$] In the annulus $A_b$, the orientation of $a$ goes from $b_1$ to $b_2$.

\item[$\bullet$] In the annulus $A_c$, the orientation of $d$ goes from $c_2$ to $c_1$.

\item[$\bullet$] The oriented boundary of the disk $M\setminus (A_b\cup A_c)$ is the concatenation
$b_1*c_1*(-b_2)*(-c_2)$, i.e. $[b_1][c_1][b_2]^{-1}[c_2]^{-1}=1$ in $\pi_1(M,W)$.

\end{itemize}

\end{definition}

\ex If $\De=\{a,b_1,b_2,c_1,c_2,d\}$ is a parallelogram decomposition for a pair $(M,\omega)$, where $\omega$ is a holomorphic $1$-form with
double zero at $W$, then $\De$ is an admissible decomposition for the pair $(M,W)$. Note that, on a fixed translation surface in $\H(2)$, there
are admissible decompositions which can not be realized as parallelogram decompositions.

\subsection{Projection to the sphere and associated symplectic homology basis}
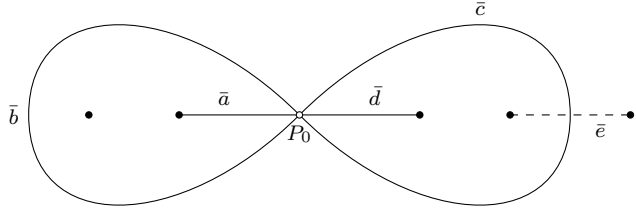
\begin{figure}[htb]
\begin{center}

\begin{tikzpicture}[scale=0.4]

\draw (0,0) .. controls (-4,4) and (-9,4) .. (-9,0); \draw (0,0) .. controls (-4,-4) and (-9,-4) .. (-9,0);

\draw (0,0) .. controls (4,4) and (9,4) .. (9,0); \draw (0,0) .. controls (4,-4) and (9,-4) .. (9,0);

\draw (-4,0) -- (0,0) -- (4,0); \draw[dashed] (7,0) -- (11,0);

\filldraw[black]  (-4,0) circle (3pt) (-7,0) circle (3pt) (4,0) circle (3pt) (7,0) circle (3pt) (11,0) circle (3pt);

\filldraw[fill=white] (0,0) circle (3pt);

\draw (0,0) node[below] {$\scriptstyle P_0$} (-9,0) node[left] {$\scriptstyle \bar{b}$} (-2.5,0) node[above] {$\scriptstyle \bar{a}$} (2.5,0) node[above] {$\scriptstyle \bar{d}$} (6,3) node[above] {$\scriptstyle \bar{c}$};

\draw (10,0) node[below] {$\scriptstyle \bar{e}$};

\end{tikzpicture}
\end{center}
\caption{Projection of an admissible decomposition to $\CP^1$.}
\label{fig:decomp-curves:CP1}

\end{figure}

Let $\rho : M \ra \CP^1$ be the two-sheeted branched  cover from $M$ onto $\CP^1$ ramified at the Weierstrass points of $M$. By definition, we have $\rho(P)=\rho(P')$ if and only if $P'\in \{P,\tau(P)\}$. Let $P_0,\dots, P_5$ denote the images of the Weierstrass points of $M$ by $\rho$, with $P_0= \rho(W)$. Let $\De=\{a,b_1,b_2,c_1,c_2,d\}$ be an admissible decomposition for the pair $(M,W)$. The projections of the curves in $\De$ satisfy (see Figure~\ref{fig:decomp-curves:CP1}).

\begin{itemize}
\item[$\bullet$] $\rho(a)=\bar{a}$ is a simple arc joining $P_0$ to another point in $\{P_1,\dots,P_5\}$. Without loss of generality, we can assume that the endpoints of $\bar{a}$ are $P_0$ and $P_1$.

\item[$\bullet$] $\rho(b_1)=\rho(b_2)=\bar{b}$ is  a simple closed curve with basepoint at $P_0$, $\bar{b}$ is the boundary of an open disc ${\rm D}_1$ which contains $\bar{a}\sm \{P_0\}$ and two points of the set $\{P_1,\dots,P_5\}$. We can assume that ${\rm D}_1\cap \{P_1,\dots,P_5\}=\{P_1,P_2\}$.

\item[$\bullet$] $\rho(d)=\bar{d}$ is a simple arc disjoint from ${\rm D}_1$ joining $P_0$ to a point in $\{P_3,\dots,P_5\}$. We can assume that the endpoints of $\bar{d}$ are $P_0$ and $P_3$.

\item[$\bullet$] $\rho(c_1)=\rho(c_2)=\bar{c}$ is a simple closed curve with basepoint $P_0$ disjoint from ${\rm D}_1$, $\bar{c}$ is the boundary of an open disc ${\rm D}_2$ which contains $\bar{d}\sm \{P_0\}$ and two points in the set $\{P_1,\dots,P_5\}$. We can assume that ${\rm D}_2\cap \{P_1,\dots,P_5\} =\{P_3,P_4\}$.
\end{itemize}

\rem The orientations of the curves in $\De$ do not determine the orientations of the curves $\{\bar{a},\bar{b},\bar{c},\bar{d}\}$.

\medskip

Let $\bar{b}_*$ be a simple arc contained in ${\rm D}_1$ which joins $P_1$ to $P_2$ such that $\bar{b}_*\cap \bar{a}=\{P_1\}$. Observe that $b=\rho^{-1}(\bar{b}_*)$ is a simple closed curve in $M$, freely homotopic to $b_1$ and $b_2$, and we have $\tau(b)=-b$. Similarly, let $\bar{c}_*$ be a simple arc in ${\rm D}_2$ joining $P_3$ to $P_4$ such that $\bar{c}_*\cap \bar{d}=\{P_3\}$, then $c=\rho^{-1}(\bar{c}_*)$ is a simple closed curve freely homotopic to $c_1,c_2$ such that $\tau(c)=-c$. 

Now, let $\bar{e}$ be a simple arc in $\CP^1$ which joins $P_5$ to $P_4$ such that $\bar{e}\cap(\bar{a}\cup\bar{b}\cup\bar{d})=\vide$ (see Figure~\ref{fig:decomp-curves:CP1}), then $e=\rho^{-1}(\bar{e})$ is a simple closed curve such that $\tau(e)=-e$. Note that $\bar{e}$ is unique up to homotopy with fixed endpoints in $\CP^1\sm (\bar{a}\cup\bar{b}\cup\bar{d})$. 

\begin{figure}[htb]
\centering
\begin{tikzpicture}[scale=0.5]
\draw (4,8) arc (90:270: 4 and 4);
\draw (2,4) .. controls (3,5) and (4,5) .. (5,4);
\draw (2,4) .. controls (3,3) and (4,3) .. (5,4);
\draw (4,8) -- (8,8) (4,0) -- (8,0);
\draw (5,4) .. controls (6,6) and (6,8) .. (5,8);
\draw (5,4) .. controls (6,2) and (6,0) .. (5,0);
\draw[dashed] (5,4) .. controls (4,6) and (4,8) .. (5,8);
\draw[dashed] (5,4) .. controls (4,2) and (4,0) .. (5,0);

\draw (8,0) arc (-90:90: 4 and 4);
\draw (5,4) .. controls (8,7) and (11,7) .. (11,4);
\draw (5,4) .. controls (8,1) and (11,1) .. (11,4);
\draw[dashed] (5,4) .. controls (8.5,7.5) and (11.5,7.5) .. (11.5,4);
\draw[dashed] (5,4) .. controls (8.5,0.5) and (11.5,0.5) .. (11.5,4);

\draw (5,4) .. controls (6,3.5) and (6,3.5) .. (7,4);
\draw[dashed] (5,4) .. controls (6,4.5) and (6,4.5) .. (7,4);

\draw (7,4) .. controls (8,5) and (9,5) .. (10,4);
\draw (7,4) .. controls (8,3) and (9,3) .. (10,4);

\draw[dashed] (0,4) .. controls (0.5,4.5) and (1.5,4.5) .. (2,4);
\draw (0,4) .. controls (0.5,3.5) and (1.5,3.5) .. (2,4);

\draw[dashed] (10,4) .. controls (10.5,4.5) and (11.5,4.5) .. (12,4);
\draw (10,4) .. controls (10.5,3.5) and (11.5,3.5) ..  (12,4);

\filldraw[fill=black] (5,4) circle (3pt);

\draw (1,3.2) node {$\scriptstyle b$};
\draw (3.5,5) node {$\scriptstyle a$};
\draw (8.5,5) node {$\scriptstyle c$};
\draw (12,4) node[right] {$\scriptstyle e$};

\end{tikzpicture}
\caption{Symplectic homology basis associated to an admissible decomposition.}
\label{fig:decomp-curv:adm-decomp}
\end{figure}
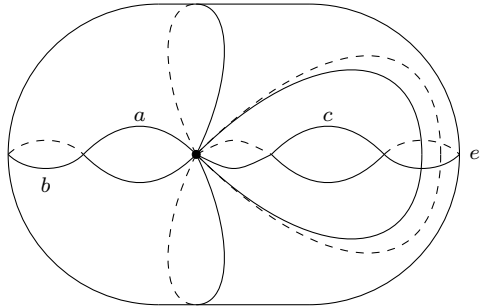

Recall that $a,d$ are already oriented, and we have an orientation for $b$ which is induced by  the orientation of $b_1$ and $b_2$. Choose the orientation of $e$ so that $\langle c, e \rangle =1$, then $(a,b,d,e)$ is a symplectic basis of $H_1(M,\Z)$. We will call $(a,b,c,e)$ the basis associated to the decomposition $\De$. It is easy to see that this basis is also the one described in Section~\ref{ThetSect}. Note also that if $\De$ is a parallelogram decomposition, then the two definitions of associated symplectic basis of $H_1(M,\Z)$ agree (see Figure~\ref{fig:decomp-curv:adm-decomp}).

\subsection{Elementary moves}

Let $\De=\{a,b_1,b_2,c_1,c_2,d\}$ be an admissible decomposition for the pair $(M,W)$. Let $(\bar{a},\bar{b},\bar{c},\bar{d},\bar{e}, {\rm D}_1,{\rm D}_2$ be as above. 

\begin{itemize}
\item[$\bullet$] \underline{The $T$ move:} \hspace{0.25cm} Let $\bar{a}'$ be a simple arc joining $P_0$ to $P_2$ such that $(\bar{a}\sm\{P_0\})\subset {\rm D}_1$, and $\bar{a}'\cap \bar{a}=\{P_0\}$. The preimage $a'$ of $\bar{a}'$ in $M$ is a simple closed curves contained in the annulus bounded by $b_1,b_2$ which satisfies $\tau(a')=-a'$. By choosing an appropriate orientation for $a'$, we see that the family $\De'=\{a',b_1,b_2,c_1,c_2,d\}$ is an admissible decomposition for $(M,W)$. The simplectic homology basis associated to $\De'$ is $(a',b,c,e)$, and we have

$$\left( \begin{array}{c} a' \\ b \\ c \\ e \end{array}\right) =\left(\begin{array}{rrrr} 1 & \pm 1 & 0 & 0 \\ 0 & 1 & 0 & 0 \\ 0 & 0 & 1  & 0 \\ 0 & 0 & 0 & 1 \end{array}\right) \cdot \left(\begin{array}{c} a \\ b \\ c \\ d \end{array} \right).$$

We call this transformation the $T$ move (see Figure~\ref{fig:dec-curv:Tmove}).

\begin{figure}[htb]
\begin{center}
\begin{tikzpicture}[scale=0.4]

\draw (0,0) .. controls (-4,4) and (-9,4) .. (-9,0); \draw (0,0) .. controls (-4,-4) and (-9,-4) .. (-9,0);

\draw (0,0) .. controls (4,4) and (9,4) .. (9,0); \draw (0,0) .. controls (4,-4) and (9,-4) .. (9,0);

\draw (0,0) .. controls (-4,2) and (-7,2) .. (-7,0);

\draw (-4,0) -- (0,0) -- (4,0); 

\filldraw[black]  (-4,0) circle (3pt) (-7,0) circle (3pt) (4,0) circle (3pt) (7,0) circle (3pt); 

\filldraw[fill=white] (0,0) circle (3pt);

\draw (0,0) node[below] {$\scriptstyle P_0$} (-9,0) node[left] {$\scriptstyle \bar{b}$} (-2.5,0) node[below] {$\scriptstyle \bar{a}$} (2.5,0) node[above] {$\scriptstyle \bar{d}$} (6,3) node[above] {$\scriptstyle \bar{c}$};


\draw (-6.8,1) node[left] {$\scriptstyle \bar{a}'$};

\end{tikzpicture}
\end{center}
\caption{$T$ move for an admissible decomposition.}
\label{fig:dec-curv:Tmove}
\end{figure}
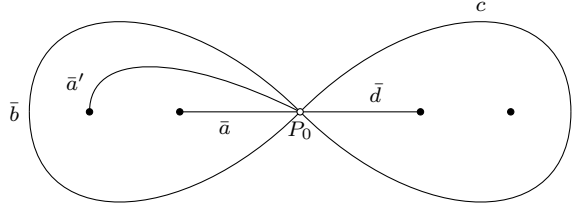.

\item[$\bullet$] \underline{The $S$ move:} We define the $S$ move for admissible decompositions in the same way as parallelogram decompositions, that is, an $S$ move transforms $\De$ into the family $\De'=\{d,-c_2,-c_1,b_1,b_2,-a\}$. Let $\bar{e}'$ be a simple arc in $\CP_1$ joining $P_5$ to $P_2$ such that $\bar{e}'\cap(\bar{a}\cup\bar{c}\cup\bar{d})=\vide$, then the symplectic homology basis associated to $\De'$ is $(d,-c,b,e')$ where $e'=\rho^{-1}(\bar{e}')$. It is easy to check that the symplectic homology bases associated to $\De$ and $\De'$ are related by the matrix $S$. 

\item[$\bullet$] \underline{The $R$ move:} \hspace{0.25cm} Let $\bar{c}'$ be a simple closed curve in $\CP^1$ with basepoint $P_0$ which is disjoint from ${\rm D}_1$ and bounds an open disc ${\rm D}'_2$ such that (see Figure~\ref{fig:dec-curv:Rmove})

\begin{itemize}
\item[.] ${\rm D}'_2\cap\{P_1,\dots,P_5\}=\{P_3,P_5\}$,

\item[.] $(\bar{d}\sm\{P_0\})\subset {\rm D}'_2$.
\end{itemize}

The preimage of $\bar{c}'$ in $M$ is the union of two simple closed curves $c'_1,c'_2$ with basepoint $W$. By choosing an appropriate orientation or $c'_1,c'_2$, we see that the family $\De'=\{a,b_1,b_2,c'_1,c'_2,d\}$ is an admissible decomposition. Let $\bar{c}'_*$ be a simple arc in ${\rm D}'_2$ joining $P_3$ to $P_5$, and $c'$ be the preimage of $\bar{c}'_*$ in $M$. Then $(a,b,c',e)$ is the symplectic homology basis associated to $\De'$. One can easily check that 

$$\left(\begin{array}{c} a \\ b\\ c'\\ e\end{array}\right) = \left(\begin{array}{rrrr} 1 & 0 & 0 & 0 \\ 0 & 1 & 0 & 0 \\ 0 & 0 & 1 & \pm 1\\ 0 & 0 & 0 & 1 \end{array}\right)\cdot \left(\begin{array}{c} a \\ b \\ c\\ e \end{array}\right).$$

\begin{figure}[htb]
\centering
\begin{tikzpicture}[scale=0.35]

\draw (0,0) .. controls (-4,4) and (-9,4) .. (-9,0); \draw (0,0) .. controls (-4,-4) and (-9,-4) .. (-9,0);

\draw (0,0) .. controls (4,4) and (9,4) .. (9,0); \draw (0,0) .. controls (4,-4) and (9,-4) .. (9,0);

\draw[red] (0,0) .. controls (3,2) and (5,1).. (6,0); \draw[red] (6,0) .. controls (7,-1) and (9,-1) .. (10,0); \draw[red] (10,0) .. controls
(11,1) and (13,1) .. (13,-1); \draw[red] (0,0) .. controls (4,-6) and (13,-6) .. (13,-1);

\draw (-4,0) -- (0,0) -- (4,0); \draw[dashed] (7,0) -- (11,0);

\filldraw[black]  (-4,0) circle (3pt) (-7,0) circle (3pt) (4,0) circle (3pt) (7,0) circle (3pt) (11,0) circle (3pt);

\filldraw[fill=white] (0,0) circle (3pt);

\draw (0,0) node[below] {$\scriptstyle P_0$} (-9,0) node[left] {$\scriptstyle \bar{b}$} (-2.5,0) node[above] {$\scriptstyle \bar{a}$} (2.5,0) node[below] {$\scriptstyle \bar{d}$} (6,3) node[above] {$\scriptstyle \bar{c}$};

\draw (8,0) node[above] {$\scriptstyle \bar{e}$};

\draw (13,-1) node[right] {$\scriptstyle \bar{c}'$};

\end{tikzpicture}
\caption{$R$ move for admissible decompositions.}
\label{fig:dec-curv:Rmove}
\end{figure}

\end{itemize}

The elementary moves $T,S,R$ defined above transform an admissible decomposition for the pair $(M,W)$ into another admissible decomposition for the same pair. Let us now introduce an elementary move which transforms an admissible decomposition for the pair $(M,W)$ into an admissible decomposition for another pair $(M,W')$.

Let $\bar{a}'$ be a simple arc in $\CP^1$ joining $P_5$ to $P_4$ such that $\bar{a}'\cap(\bar{a}\cup\bar{b}\cup\bar{d})=\vide$. Let $\bar{d}'$ be a simple arc joining $P_5$ to $P_2$ such that $\bar{d}'\cap(\bar{a}\cup\bar{c}\cup\bar{d})=\vide$. Let $\bar{b}'$ and $\bar{c}'$ be two simple closed curves with basepoint $P_5$ which are respectively the boundary of the open discs ${\rm D}'_1$ and ${\rm D}'_2$ satisfying  (see Figure~\ref{fig:dec-curv:Umove})

\begin{itemize}
\item[.] ${\rm D}'_1\cap {\rm D}'_2=\vide, {\rm D'}_i\cap {\rm D}_i=\vide, i=1,2$,

\item[.] ${\rm D}'_1\cap\{P_0,\dots,P_4\}=\{P_3,P_4\}, {\rm D}'_2\cap\{P_0,\dots,P_4\}=\{P_1,P_2\}$,

\item[.] $\inter(\bar{a}')\subset {\rm D}'_1, \inter(\bar{d}')\subset {\rm D}'_2.$

\end{itemize}

Then the preimages of $\bar{a}',\bar{b}',\bar{c}',\bar{d}'$ is an admissible decomposition $\De'$ for the pair $(M,\rho^{-1}(P_5))$. Actually, the family $(\bar{a}',\bar{b}',\bar{c},\bar{d}')$ gives rise to two admissible decompositions, each of which is determined by the orientation of $a'=\rho^{-1}(\bar{a}')$. The symplectic homology bases associated  those decompositions are related by $-\Id$. Let $(a',b',c',e')$ be the symplectic homology basis associated to $\De'$,  we have

$$\left(%
\begin{array}{c}
  a' \\ b' \\ c' \\ e' \\
\end{array}%
\right) = \pm\left(%
\begin{array}{cccc}
  0 & 0 & 0 & -1 \\
  0 & 0 & 1 & 0 \\
  0 & -1 & 0 & 0 \\
  1 & 0 & 0 & 0\\
\end{array}%
\right)\cdot \left(%
\begin{array}{c}
  a \\ b \\ c \\ e \\
\end{array}%
\right)=  U^{\pm 1} \cdot \left(%
\begin{array}{c}
  a \\ b \\ c \\ e \\
\end{array}%
\right).$$

\noindent We will call the transformation from $\De$ into $\De'$ the $U$ move.

\begin{figure}[htb]
\begin{minipage}[t]{0.45\linewidth}
\centering
\begin{tikzpicture}[scale=0.4]
\draw (0,0) .. controls (-4,4) and (-9,4) .. (-9,0); \draw (0,0) .. controls (-4,-4) and (-9,-4) .. (-9,0);

\draw (0,0) .. controls (4,4) and (9,4) .. (9,0); \draw (0,0) .. controls (4,-4) and (9,-4) .. (9,0);

\draw (-4,0) -- (0,0) -- (4,0); 

\filldraw[black]  (-4,0) circle (3pt) (-7,0) circle (3pt) (4,0) circle (3pt) (7,0) circle (3pt);

\filldraw[fill=white] (0,0) circle (3pt);

\draw (0,0) node[above] {$\scriptstyle P_0$} (-9,0) node[left] {$\scriptstyle  \bar{b}$} (-1.5,0) node[above] {$\scriptstyle  \bar{a}$} (1.5,0) node[above] {$\scriptstyle \bar{d}$} (6,3) node[above] {$\scriptstyle \bar{c}$};


\draw[red, dashed] (0,-5) .. controls (-1,0) and (-4,2) .. (-6,2); \draw[red, dashed] (-6,2) .. controls (-7.5,2) and (-8,1) .. (-8,0);
\draw[red, dashed] (-8,0) .. controls (-8,-1) and (-7,-4) .. (0,-5);

\draw[red, dashed] (0,-5) .. controls (1,0) and (4,2) .. (6,2); \draw[red, dashed] (6,2) .. controls (7.5,2) and (8,1) .. (8,0); \draw[red,
dashed] (8,0) .. controls (8,-1) and (7,-4) .. (0,-5);

\draw[red, dashed] (-7,0) -- (0,-5) -- (7,0);

\filldraw (0,-5) circle (3pt);

\draw (0,-5) node[below] {$\scriptstyle P_5$}; \draw[red] (-4,-4) node[below] {$\scriptstyle \bar{c}'$} (4,-4) node[below] {$\scriptstyle \bar{b}'$}; \draw[red] (-5,-1.5) node[above right] {$\scriptstyle \bar{d}'$} (5,-1.5) node[above left] {$\scriptstyle \bar{a}'$};

\draw (-4,0) node[above] {$\scriptstyle P_1$} (-7,0) node[above] {$\scriptstyle P_2$} (4,0) node[above] {$\scriptstyle P_3$} (7,0) node[above] {$\scriptstyle P_4$};

\end{tikzpicture}
\end{minipage}
\begin{minipage}[t]{0.45\linewidth}
\centering
\begin{tikzpicture}[scale=0.5]
\draw (4,8) arc (90:270: 4 and 4);
\draw (4,8) -- (8,8) (4,0) -- (8,0);
\draw (8,0) arc (-90:90: 4 and 4);

\draw (3.5,4) ellipse (1.5 and 1) (8.5,4) ellipse (1.5 and 1);

\draw[dashed] (0,4) .. controls (0.5,4.5) and (1.5,4.5) .. (2,4);
\draw (0,4) .. controls (0.5,3.5) and (1.5,3.5) .. (2,4);

\draw[dashed] (10,4) .. controls (10.5,4.5) and (11.5,4.5) .. (12,4);
\draw (10,4) .. controls (10.5,3.5) and (11.5,3.5) ..  (12,4);

\draw (5,4) .. controls (6,3.5) and (6,3.5) .. (7,4);
\draw[dashed] (5,4) .. controls (6,4.5) and (6,4.5) .. (7,4);

\draw (8,0.3) arc (-90:90: 4 and 3.7);
\draw (8,0.3) -- (4.5,0.3) (8,7.7) -- (4.5,7.7);

\draw (4.5,0.3) .. controls (4,0.3) and (4,3) .. (3.5,3);
\draw[dashed] (4.5,0) .. controls (4,0) and (3,3) .. (3.5,3);
\draw (4.5,7.7) .. controls (4,7.7) and (4,5) .. (3.5,5);
\draw[dashed] (3.5,5) .. controls (3,5) and (4,8) .. (4.5,8);

\draw (12,4) .. controls (10,8) and (6,8) .. (6,4);
\draw (12,4) .. controls (10,0) and (6,0) .. (6,4);

\draw[dashed] (12,4) .. controls (9,8) and (5.5,8) .. (5.5,4);
\draw[dashed] (12,4) .. controls (9,0) and (5.5,0) .. (5.5,4);

\draw (1,3.2) node {$\scriptstyle c'$};
\draw (3.5,5) node[below] {$\scriptstyle d'$};
\draw (8.5,5) node[above] {$\scriptstyle b'$};
\draw (11,3.3) node[left] {$\scriptstyle a'$};

\filldraw[fill=white] (5,4) circle (3pt);
\filldraw[fill=black] (12,4) circle (3pt);

\end{tikzpicture}
\end{minipage}
\caption{$U$ move for admissible decompositions.}
\label{fig:dec-curv:Umove}
\end{figure}

From the definitions, the following lemma is clear

\begin{lemma}\label{lm:Gam:Adm}
Let $\De$ be an admissible decomposition for the pair $(M,W)$, and $\gamma$ an element of the group $\Gamma$. Then there exists an admissible decomposition $\De'$ for the same pair such that the symplectic homology bases associated to $\De$ and $\De'$ are related by $\gamma$.
\end{lemma}

\medskip

From Lemma~\ref{SpGenLm}, we know that the family $\{T,S,R,U\}$ generates the group $\Sp(4,\Z)$, therefore we have

\begin{lemma}\label{lm:Sp:Adm}
Let $\De$ be an admissible decomposition for the pair $(M,W)$, and $A$ an element of $\Sp(4,\Z)$. Then there exist a Weierstrass point $W'$ of $M$, and an admissible decomposition $\De'$ for the pair $(M,W')$ such that the symplectic homology bases associated to $\De$ and $\De'$ are related by $A$.
\end{lemma}

\section{Symplectic homology bases associated to admissible decompositions and the group $\Gamma$.}\label{sect:Gam:SympBases}

Let $(M,W)$ be an element of the space $\H(2)/\C^*$. Our aim in this section is to prove the following theorem, which is the key ingredient of the proof of Theorem~\ref{th0}. 

\begin{theorem}\label{th1}

Let $\De$ and $\De'$ be two admissible decompositions for the pair $(M,W)$ with the associated symplectic homology bases $(a,b,c,e)$, and
$(a',b',c',e')$ respectively. Then there exists an element $\gamma$ in $\Gamma$ such that

$$\left(%
\begin{array}{c}
  a' \\
  b' \\
  c' \\
  e' \\
\end{array}%
\right)= \gamma \cdot \left(%
\begin{array}{c}
  a \\
  b \\
  c \\
  e \\
\end{array}%
\right).$$

\end{theorem}

In what follows, we denote by $\tau$ the hyperelliptic involution of $M$, by $\rho$ the two-sheeted branched cover from $M$ onto $\CP^1$, and by $P_0,\dots,P_5$ the images of the Weierstrass points of $M$ by $\rho$, where $P_0=\rho(W)$. We also equip $M$ with the hyperbolic metric in the conformal class of the Riemann surface structure. Note that $\tau$ is now an isometry of $M$. Recall that, for any closed curve $\alpha$ in $M$ which contains $W$, we denote by $[\alpha]$ the homotopy class of $\alpha$ in $\pi_1(M,W)$.

\subsection{Admissible decompositions with common subfamily}
Let us first prove the following

\begin{proposition}\label{PrfTh1Prop0}
Let $\De=\{a,b_1,b_2,c_1,c_2,d\}$ and $\De'=\{a',b'_1,b'_2,c'_1,c'_2,d'\}$ be two admissible decompositions for the pair $(M,W)$. Assume that  $\bar{b}=\rho(b_1)$, and $\bar{b}'=\rho(b'_1)$ are homotopic in $\pi_1(\CP^1\sm\{P_1,\dots,P_5\}, P_0)$. Then there exists an element $\gamma\in \Gamma$ such that the symplectic bases of $H_1(M,\Z)$ associated to $\De$ and $\De'$ are related by $\gamma$.

%

\end{proposition}

\dem Since $\bar{b}$ and $\bar{b}'$ are homotopic in $\pi_1(\CP^1\sm\{P_1,\dots,P_5\},P_0)$, there exists a homeomorphism $\Phi$ of $\CP^1$  isotopic to the identity relative to $\{P_0,\dots,P_5\}$ such that $\Phi(\bar{b})=\bar{b}'$. The homotopy from the identity of $\CP^1$ to $\Phi$ can be lifted to a homotopy of $M$ which is identity on the set of Weierstrass points of $M$. Therefore, we can assume that $b_1\cup b'_1=b'_1\cup b'_2$ as subsets of $M$. A priori, the orientations of $(b_1,b_2)$ and $(b'_1,b'_2)$ may not be the same, but since $-\Id \in \Gamma$, we can assume that they have the same orientation which means that $b_1=b'_1$ in $H_1(M,\Z)$. Note that the orientation of $(b_1,b_2)$ determines the orientation of $a$ by the condition $\langle a, b_1 \rangle =1$, and consequently, we get a unique numbering of the pair $(b_1,b_2)$ by the condition that $a$ goes from $b_1$ to $b_2$. 

\noindent By definition, $b_1\cup b_2$ is the boundary of an open annulus $A_b$ which contains both $a\sm\{W\}$ and $a'\sm\{W\}$, therefore  there exists an integer $n$ such that $a'$ is homotopic to the image of $a$ by $n$ Dehn twists in $A_b$. Thus, by applying the $T$ move $n$ times, we can assume that $a'=a$ as subsets of $M$. Let $(a,b,c,e)$ and $(a',b',c',e')$ be the symplectic bases of $H_1(M,\Z)$ associated to $\De$ and $\De'$ respectively.    It follows that we have the following equality in $H_1(M,\Z)$

$$\left(%
\begin{array}{c}
  a' \\
  b' \\
  c' \\
  e' \\
\end{array}%
\right)= \left(%
\begin{array}{cc}
  \Id_2 & 0 \\
  X & Y \\
\end{array}%
\right) \cdot \left(%
\begin{array}{c}
  a \\
  b \\
  c \\
  e \\
\end{array}%
\right),$$

\noindent with $X,Y \in \mathbf{M}_{2\times 2}(\Z)$. Since $\left(%
\begin{array}{cc}
  \Id_2 & 0 \\
  X & Y \\
\end{array}%
\right)$ belongs to $\Sp(4,\Z)$, simple computations show that we must have $X=0$, and $Y\in \SL(2,\Z)$. Since the group $\DS{\left(%
\begin{array}{cc}
  \Id_2 & 0 \\
  0 & \mathrm{SL}(2,\Z) \\
\end{array}%
\right)}$ is contained in $\Gamma$ (cf. Lemma~\ref{GamPropLm}), the proposition follows. \carre


\subsection{Standard decomposition}

Let us now prove the following

\begin{lemma}\label{PrfTh1Lm1}
Let $\De=\{a,b_1,b_2,c_1,c_2,d\}$ be an admissible decomposition of the pair $(M,W)$. Let $b_0$ and $c_0$  be the simple closed geodesic in the
free homotopy class of $b_1$ and $c_1$ respectively, then we have

\begin{itemize}

\item[a)] $b_0\cap c_0=\vide$.

\item[b)] $\tau(b_0)=-b_0, \tau(c_0)=-c_0$.

\item[c)] $W\notin b_0\cup c_0$.

\end{itemize}

\end{lemma}

\dem

\begin{itemize}

\item[a)] Let  $b$ (resp. $c$) be a simple closed curve freely homotopic to $b_1$ (resp. to $c_1$) which is contained in the annulus bounded by $b_1$ and $b_2$ (resp. $c_1$ and $c_2$). By construction $b$ and $c$ are freely homotopic to $b_0$ and $c_0$ respectively, and  $b\cap c=\vide$.
It is well known that $\card\{b_0\cap c_0\}$ is the intersection number $\iota(b,c)$ of the free homopoty classes of $b$ and $c$, thus we have
that $\card\{b_0\cap c_0\}=0$.

\item[b)] By definition, we see that $b_1$ is freely homotopic to $b_2$, and  $\tau(b_1)=-b_2$, therefore $\tau(b_1)$ is freely
homotopic to $-b_1$. Now, since $b_0$ is the unique simple closed geodesic in the free homotopy class of $b_1$,  $-b_0$ is then the unique
simple closed geodesic in the free homotopy class of $-b_1$. Since $\tau$ is an isometry, $\tau(b_0)$ must be simple closed geodesic in the free homotopy class of $\tau(b_1)$, hence $\tau(b_0)$ is freely homotopic to $-b_1$, therefore we must have $\tau(b_0)=-b_0$.

\item[c)] Suppose that $W\in b_0$. Since $b_0$ is freely homotopic to $b_1$, there exists $[h]\in \pi_1(M,W)$ such that $[b_0]=[h][b_1][h]^{-1}$. Note that we have $[\tau(b_0)]=[b_0]^{-1}$, hence,

$$ [h][b_1]^{-1}[h]^{-1}=[\tau(h)][\tau(b_1)][\tau(h)]^{-1}= [\tau(h)][a]^{-1}[b_1]^{-1}[a][\tau(h)]^{-1}.$$

\noindent It follows

\begin{equation}\label{HomEq1}
[a][\tau(h)]^{-1}[h][b_1]^{-1}=[b_1]^{-1}[a][\tau(h)]^{-1}[h]
\end{equation}

\noindent We deduce that $[b_1]^{-1}$ and $[a][\tau(h)]^{-1}[h]$ commute. But $[b_1]$ is a simple closed non-separating curve, therefore we have

\begin{equation}\label{HomEq2}
[a][\tau(h)]^{-1}[h]=[b_1]^{n} \mbox{ with $n\in \Z$}
\end{equation}

\noindent Recall that $\tau$ acts like $-\Id$ on $H_1(M,\Z)$, thus (\ref{HomEq2}) implies the following equality in $H_1(M,\Z)$

\begin{equation}\label{HomEq3}
nb_1-a=2h
\end{equation}

\noindent It follows that $\langle a, b_1\rangle =0 \mod 2$,  but by construction we know that $\langle a, b_1\rangle =1$, and we get a contradiction.  We can then conclude that $b_0$ does not contain $W$. The same arguments apply to $c_0$ and the lemma follows.

\end{itemize}

\carre

\rem  If a simple closed curve $g$ satisfies $\tau(g)=-g$, then $g$ contains exactly two  fixed points of $\tau$, which are Weierstrass points of
$M$.

\bigskip

Let $(g_1,g_2)$ be a pair of disjoint simple closed geodesics verifying the following property

$$(\mathscr{P}) \; \left\{%
\begin{array}{ll}
    W \notin g_i,  \\
    \tau(g_i)=-g_i, \\
\end{array}%
\right.$$

\noindent for $i=1,2$. We construct an admissible decomposition of $(M,W)$ associated to $(g_1,g_2)$ as follows:

\begin{itemize}

\item[-] Cut open $M$ along $g_1$ and $g_2$, we obtain a $4$-holed sphere $N$ which is equipped with a hyperbolic metric with geodesic boundary.
Let $g_i^+,g_i^-$ denote the boundary components of $N$ corresponding to $g_i, \; i=1,2$. The orientation of $g_i$ implies an orientation for $g^{\pm}_i$. We choose the notations  so that the orientation of $g^+_i$ agrees with the orientation induced by the orientation of $N$. Note that the hyperelliptic involution $\tau$ of $M$ induces an isometric involution $\tau'$ of $N$ which interchanges $g^+_i$ and $g_i^-$.

\item[-] Let $s_i^+$ be the shortest path in $N$ from $W$ to $g^+_i$, and  let $s_i^-$ denote $\tau'(s^+_i), \; i=1,2$. Note that the action of
$\tau$ on the tangent space at $W$ is $-\Id$, therefore $s_i^+\cup s_i^-$ is a simple geodesic arc joining $g^+_i$ to $g^-_i$. By a slight abuse of notations, we also denote by $s^+_i, s^-_i$ be the geodesic arcs in $M$ corresponding to $s^+_i$ and $s^-_i$. Let $P_i^+,P_i^-$ denote the endpoints of $s^+_i$, and $s^-_i$ in $g_i$ respectively.  Let $r_i$ be a simple arc in $g_i$ with endpoints $P_i^+,P_i^-$, then $h_i=r_i\cup s^+_i\cup s^-_i$ is a simple closed curve containing $W$. 

\item[-] We choose the orientation of $h_1$ so that $\langle h_1,g_1\rangle=1$. Observer the induced orientation for $s^+_1$  is from $W$ to $P_1^+$.  We can then find a simple closed curve $\tilde{b}$ homotopic to $s^+_1*g_1*(-s^+_1)$ in $\pi_1(M,W)$ such that

\begin{itemize}


\item[.] $\tilde{b}\cap\tau(\tilde{b})=\{W\}$,

\item[.] $\tilde{b}$ and $\tau(\tilde{b})$ bound an annulus containing $g_1$ as a waist curve.

\end{itemize}

\item[-] We choose the orientation of  $t_2$ so that $\langle g_2,h_2\rangle =1$. By assumption on the orientation of $g^{\pm}_2$, we see that the
induced orientation of $s^+_2$ is from $P^+_2$ to $W$. Let $\tilde{c}$ be a simple closed curve homotopic to $(-s^+_2)*g_2*s^+_2$ in
$\pi_1(M,W)$ such that $\tau(\tilde{c})\cap \tilde{c}=\{W\}$. Following the orientation of $g_2$, we have two cases:

\begin{enumerate}

\item $[\tilde{b}][\tilde{c}][\tau(\tilde{b})][\tau(\tilde{c})]=1$ in $\pi_1(M,W)$:  in this case, take $a=h_1, b_1=\tilde{b},
b_2=-\tau(\tilde{b}), c_1=\tilde{c}, c_2=-\tau(\tilde{c}), d=h_2$, then $\De=\{a,b_1,b_2,c_1,c_2,d\}$ is an admissible decomposition for
$(M,W)$ (see Figure~\ref{fig:standDec:case1}).

\begin{figure}[!h]
\begin{center}
  \includegraphics[width=5.5cm, clip=, bbllx=100, bblly=50, bburx=1400, bbury=1330]{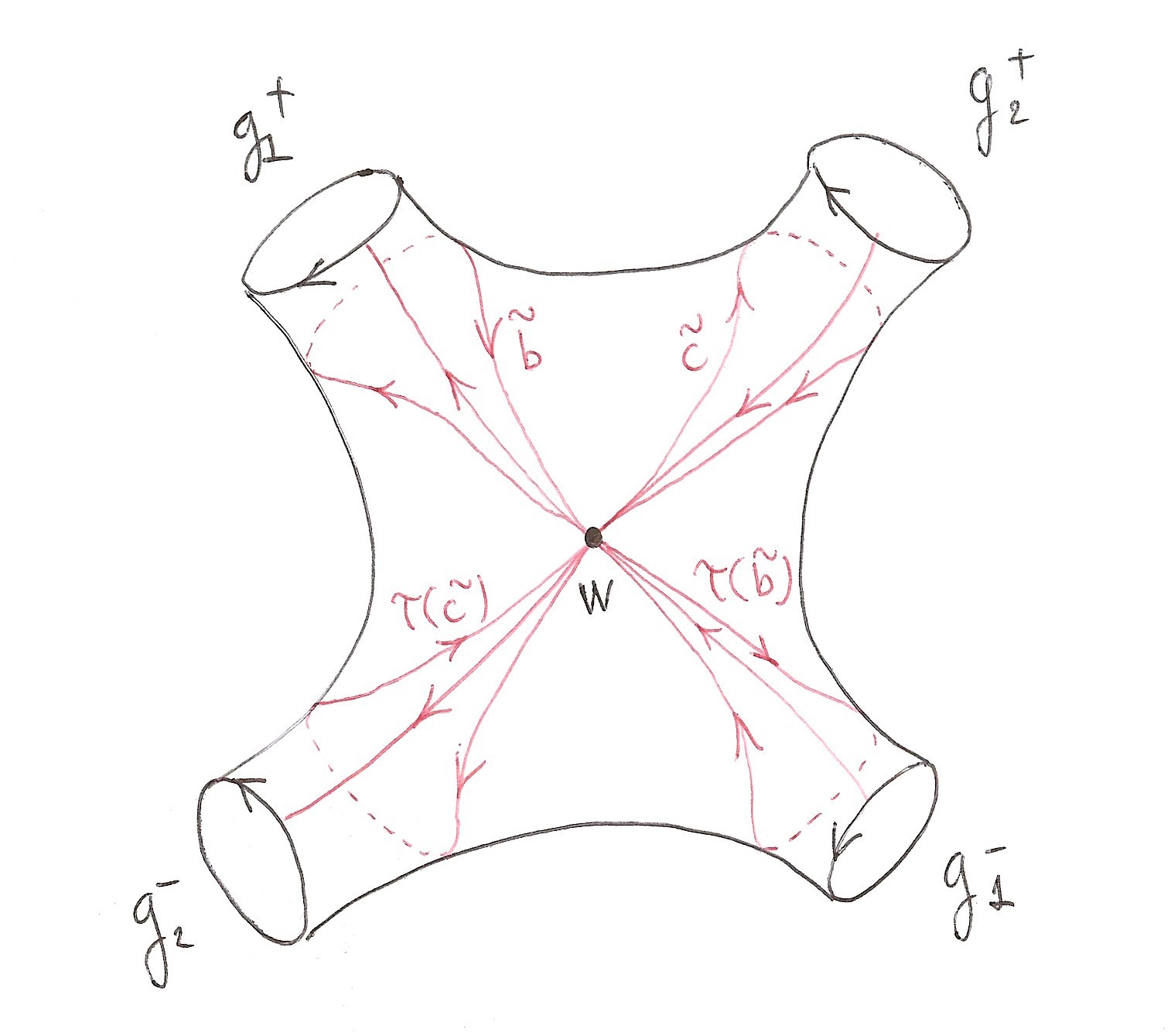}
\end{center}
\caption{Admissible decomposition associated to a pair of geodesics satisfying $(\mathscr{P})$: case 1.}
\label{fig:standDec:case1}
\end{figure}

\item $[\tilde{b}][\tau(\tilde{c})][\tau(\tilde{b})][c]=1$ in $\pi_1(M,W)$: in this case,  let $u^+_2$ be a simple arc going from
$P^+_2$ to $W$ such that $(-s^+_2)*u^+_2$ is homotopic to $\tilde{b}$ in $\pi_1(M,W)$, and $u^+_2\cap \tau(u^+_2)=\{W\}$. Then there
exists a simple closed curve $\tilde{c}'$ homotopic to $(-u^+_2)*g_2*u^+_2$ in $\pi_1(M,W)$ such that

\begin{itemize}

\item[$\bullet$] $\tilde{c}'\cap \tau(\tilde{c}')=\{W\}$,

\item[$\bullet$] $\tilde{c}'$ and $\tau(\tilde{c}')$ bound an annulus containing $g_2$ as a waist curve.
\end{itemize}

\begin{figure}[!h]
\begin{center}
  \includegraphics[width=6cm, clip=, bbllx=0, bblly=200, bburx=1400, bbury=1350]{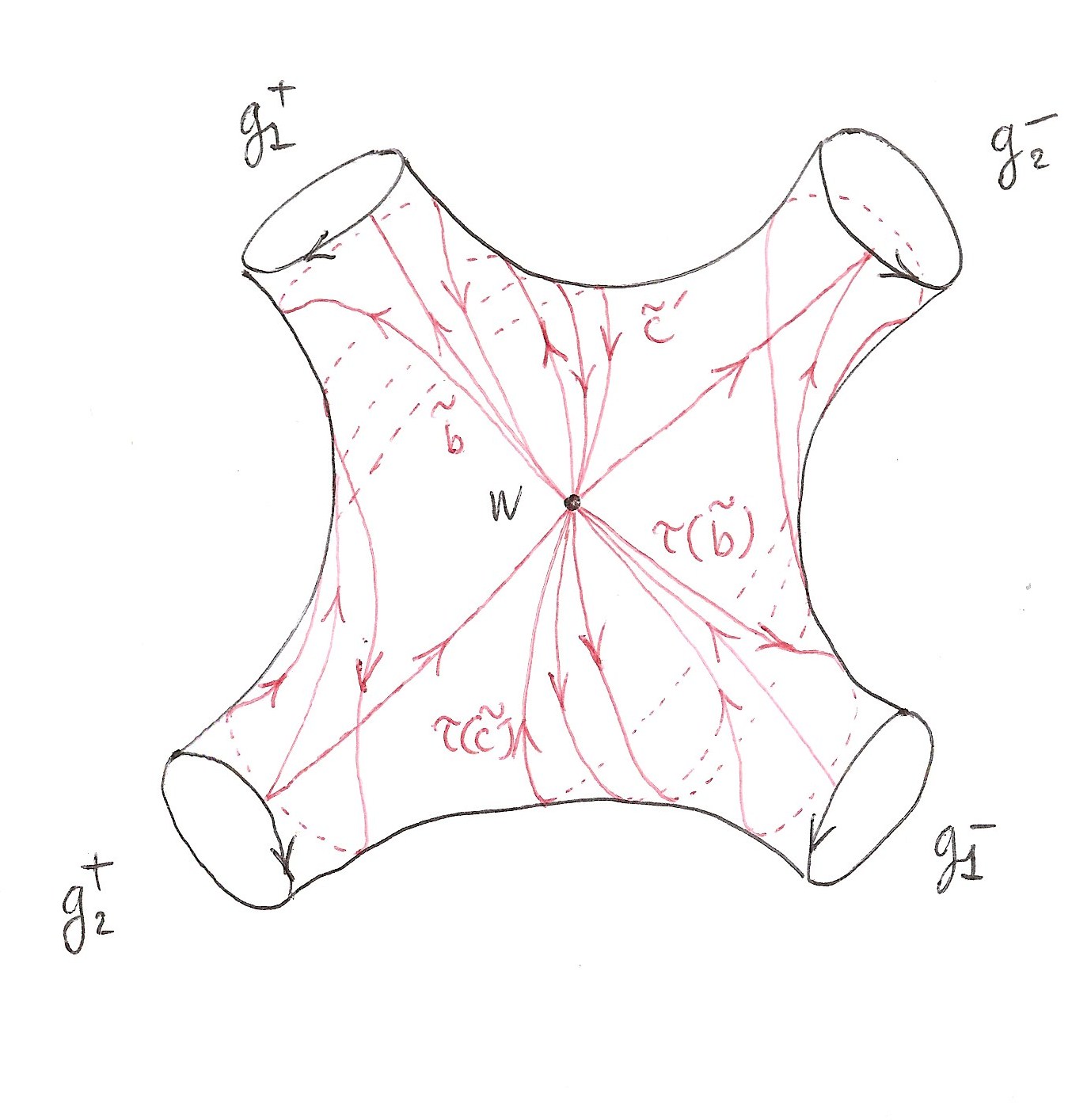}

\end{center}
\caption{Admissible decomposition associated to a pair of geodesics satisfying $(\mathscr{P})$: case 2.}
\label{fig:standDec:case2}
\end{figure}

We have

$$\begin{array}{rcl}
[\tilde{c}'] & = & [(-u_2^+)*g_2*u_2^+],\\
             & = & [(-u_2^+)*s_2^+][(-s_2^+)*g_2*s_2^+][(-s_2^+)*u_2^+],\\
             & = & [\tilde{b}]^{-1}[\tilde{c}][\tilde{b}].

\end{array}$$

\noindent Therefore,

$$[\tilde{b}][\tilde{c}'][\tau(\tilde{b})][\tau(\tilde{c}')]=1.$$

\noindent Let $h'_2$ denote the simple closed curve $r_2\cup \tau(u^+_2)\cup u^+_2$, we choose the orientation of $h'_2$ so that $\langle
g_2,h'_2\rangle=1$. It follows that, if we take $a=h_1, b_1=\tilde{b}, b_2=-\tau(\tilde{b}), c_1=\tilde{c}', c_2=-\tau(\tilde{c}')$, and
$d=h'_2$, then $\De=\{a,b_1,b_2,c_1,c_2,d\}$ is an admissible decomposition for $(M,W)$ (see Figure~\ref{fig:standDec:case2}).

\end{enumerate}

In both cases, we will call $\De$ a standard decomposition associated to the pair of (oriented) geodesics $(g_1,g_2)$.

\end{itemize}

From Lemma \ref{PrfTh1Lm1}, we know that if $\De=\{a,b_1,b_2,c_1,c_2, d\}$ is an admissible decomposition for the pair $(M,W)$, and $b_0$ (resp.
$c_0$) is the simple closed geodesic in the free homotopy class of $b_1$ (resp. $c_1$), then the pair $(b_0,c_0)$ satisfies Property
$(\mathscr{P})$. Hence, we can consider the standard decompositions associated to $(b_0,c_0)$. The following proposition tells us that the
symplectic bases of $H_1(M,\Z)$ associated to the two decompositions are related by an element of the group $\Gamma$.

\begin{proposition}\label{StdDecpProp}

Let $\De=\{a,b_1,b_2,c_1,c_2,d\}$ be an admissible decomposition for the pair $(M,W)$. Let $b_0, c_0$ be the simple closed geodesics in the free
homotopy classes of $b_1$ and $c_1$ respectively. Let $\hat{\De}=\{\hat{a},\hat{b}_1,\hat{b}_2,\hat{c}_1,\hat{c}_2,\hat{d}\}$ be a standard
decomposition associated to the pair $(b_0,c_0)$, then the symplectic bases of $H_1(M,\Z)$ associated to $\De$ and $\hat{\De}$ are related by an
element in $\Gamma$.

\end{proposition}

\dem Since $b_1$ and $\hat{b}_1$ are freely homotopic, there exists a closed curve $h$ with base point $W$ such that $[\hat{b}_1]=[h]^{-1}[b_1][h]$ in $\pi_1(M,W)$. We will show that $h\in \Z a \oplus \Z b_1 \oplus \Z c_1 $ in $H_1(M,\Z)$. By definition, we have

\begin{eqnarray*}
  [\tau(\hat{b}_1)] & = & [\hat{a}]^{-1}[\hat{b}_1]^{-1}[\hat{a}]\\
 \Rightarrow \:  [\tau(h)]^{-1}[\tau(b_1)][\tau(h)] & = & [\hat{a}]^{-1}[h]^{-1}[b_1]^{-1}[h][\hat{a}]\\
 \Rightarrow \:  [\tau(h)]^{-1}[a]^{-1}[b_1]^{-1}[a][\tau(h)] & = & [\hat{a}]^{-1}[h]^{-1}[b_1]^{-1}[h][\hat{a}]\\
 \Rightarrow  \: [b_1]^{-1}([a][\tau(h)][\hat{a}]^{-1}[h]^{-1}) & = & ([a][\tau(h)][\hat{a}]^{-1}[h]^{-1})[b_1]^{-1}\\
\end{eqnarray*}

\noindent It follows that $[b_1]$ and $[a][\tau(h)][\hat{a}]^{-1}[h]^{-1}$ commute. Since $b_1$ is a simple closed curve, there exists $k\in \Z$
such that

$$[a][\tau(h)][\hat{a}]^{-1}[h]^{-1}= [b_1]^k.$$

\noindent Therefore, in $H_1(M,\Z)$, we have $ \hat{a}=a-kb_1-2h$. We know that $\langle\hat{a},\hat{c}_1\rangle=\langle\hat{a},c_1\rangle=0$, and $\langle a,c_1\rangle=\langle b_1,c_1\rangle=0$, hence, $\langle h,c_1\rangle=0$,  which implies that

$$h\in c_1^{\bot}= \Z a \oplus \Z b_1 \oplus \Z c_1.$$

Let $(a,b,c,e)$ and $(\hat{a},\hat{b},\hat{c},\hat{e})$ be the symplectic bases of $H_1(M,\Z)$ associated to $\De$ and $\hat{\De}$ respectively.
We know that there exists $\gamma \in \mathrm{Sp}(4,\Z)$ such that

$$\left(%
\begin{array}{c}
  \hat{a} \\ \hat{b} \\ \hat{c} \\ \hat{e} \\
\end{array}%
\right) = \gamma \cdot \left(%
\begin{array}{c}
  a \\   b \\ c \\ e \\
\end{array}%
\right).$$

\noindent Since in $H_1(M,\Z)$ we have $\hat{b}=b, \hat{c}=c$, and $\hat{a}=a+kb+2h \in \Z a \oplus \Z b \oplus \Z c$, it follows that $\gamma$
is of the form

$$\gamma=\left(%
\begin{array}{cccc}
  x & y & z & 0 \\
  0 & 1 & 0 & 0 \\
  0 & 0 & 1 & 0 \\
  x' & y' & z' & t' \\
\end{array}%
\right).$$

\noindent Now, $\langle \hat{a},\hat{b}\rangle =1$ implies  $x=1$, and since $z$ is the $c$-coordinate of $2h$ in $\Z a \oplus \Z b \oplus \Z
c$, we have

$$ \hat{a}= a +mb +2\ell c, \: \text{ with } m,\ell \in \Z.$$

\noindent It follows,

$$ \begin{array}{ccl}
\langle \hat{b},\hat{e}\rangle= 0 & \Rightarrow & x'=0,\\
\langle \hat{c},\hat{e}\rangle= 1 & \Rightarrow & t'=1,\\
\langle \hat{a},\hat{e}\rangle = 0 & \Rightarrow & y'=-z=-2\ell.\\
\end{array}$$

\noindent We deduce that

$$ \gamma=\left(%
\begin{array}{cccc}
  1 & m & 2\ell & 0 \\
  0 & 1 & 0 & 0 \\
  0 & 0 & 1 & 0 \\
  0 & -2\ell & n & 1 \\
\end{array}%
\right),$$

\noindent with $\ell, m, n$ in $\Z$. The proposition follows from Lemma \ref{StdDecpLm} here below. \carre

\begin{lemma}\label{StdDecpLm}
For any integers $\ell, m, n$, the matrix $\gamma = \left(%
\begin{array}{cccc}
  1 & m & 2\ell & 0 \\
  0 & 1 & 0 & 0 \\
  0 & 0 & 1 & 0 \\
  0 & -2\ell & n & 1 \\
\end{array}%
\right)$ belongs to the group $\Gamma$.\\

\end{lemma}

\dem We have shown that the group $\left(%
\begin{array}{cc}
  \Id_2 & 0\\
  0 & \mathrm{SL}(2,\R) \\
\end{array}%
\right)$ is included in $\Gamma$. Thus

$$X=S\cdot \left(%
\begin{array}{cc}
  \Id_2 & 0 \\   0 & -\Id_2 \\
\end{array}%
\right) \cdot S \cdot \left(%
\begin{array}{cc}
  \Id_2 & 0 \\   0 & -\Id_2 \\
\end{array}%
\right)=\left(%
\begin{array}{cccc}
  1 & 0 & 2 & 0 \\
  0 & 1 & 0 & 0 \\
  0 & 0 & 1 & 0 \\
  0 & -2 & 0 & 1 \\
\end{array}%
\right) \in \Gamma.$$

Set

$$Y=S^{-1}T^{-1}S=\left(%
\begin{array}{cccc}
  1 & 0 & 0 & 0 \\
  0 & 1 & 0 & 0 \\
  0 & 0 & 1 & 0 \\
  0 & 0 & 1 & 1 \\
\end{array}%
\right).$$

\noindent Now, straight computations show that $T,X$, and $Y$ commute, and $\gamma=T^{m}X^\ell Y^n$. \carre

\begin{corollary}\label{StdDecpCor}
If $\De=\{a,b_1,b_2,c_1,c_2,d\}$, and $\De'=\{a',b'_1,b'_2,c'_1,c'_2,d'\}$ are two admissible decomposition of the pair $(M,W)$ such that $b_1$
is freely homotopic to $b'_1$, and $c_1$ is freely homotopic to $\pm c'_1$. Then the symplectic  bases of $H_1(M,\Z)$ associated to $\De$ and $\De'$ are related by an element in $\Gamma$.
\end{corollary}

\dem Let $b_0,c_0,c'_0$ be the simple closed geodesics in the free homotpoty classes of $b_1,c_1$, and $c'_1$ respectively. Let $\hat{\De}=\{\hat{a},\hat{b}_1,\hat{b}_2,\hat{c}_1,\hat{c}_2,\hat{d}\}$, and $\hat{\De}'=\{\hat{a}'
,\hat{b}'_1,\hat{b}'_2,\hat{c}'_1, \hat{c}'_2,\hat{d}'\}$ be some standard decompositions associated to the pairs $(b_0,c_0)$ and $(b_0,c'_0)$
respectively.  From Proposition \ref{StdDecpProp}, the symplectic bases of $H_1(M,\Z)$ associated to $\De$ and $\hat{\De}$ (resp. $\De'$ and
$\hat{\De}'$) are related by an element of $\Gamma$.

\noindent By assumption, we have $c_0=\pm c'_0$. From the construction of standard decompositions, we can assume that $\hat{b}_1=\hat{b}'_1, \hat{b}_2=\hat{b}'_2$, and the corollary follows from Proposition \ref{PrfTh1Prop0}.\carre

\subsection{Reducing the number of intersections}
By Proposition \ref{StdDecpProp}, we can now restrict ourselves into the case of standard decompositions. We start by proving the following key
lemma

\begin{lemma}\label{Key1Lm}
Let $g,g_1,g_2$ be three simple closed geodesics of $M$ verifying Property $(\P)$. Assume that

$$\left\{%
\begin{array}{l}
    g_1\cap g_2=\vide, \\
    \card\{g \cap (g_1\cup g_2)\}=n >1.  \\
\end{array}%
\right.$$

\noindent Then there exists a simple closed geodesic $g_3$ verifying Property $(\P)$ such that $\card\{g_3\cap (g_1\cup g_2)\}=1$, and suppose that $g_1\cap g_3=\vide$, then $\card\{g\cap(g_1\cup g_3)\}<\card\{g\cap(g_1\cup g_2)\}$. Moreover, for $i=1,2$, if $g\cap g_i=\vide$ then we can find $g_3$ such that $g_3\cap g_i=\vide$. 
\end{lemma}

\dem We know that each of the curves $g_1,g_2$ contains two Weierstrass points. Let $W'$ be the other Weierstrass point of $M$ which is not
contained in $g_1\cup g_2$. We have two possibilities:

\medskip

\noindent \underline{\bf Case 1: $W'\in g$.}  Let $s$ be the segment of $g$ which contains $W'$ with endpoints in $g\cap (g_1\cup g_2)$. We
denote by $Q_1,Q_2$ the two endpoints of $s$, and choose the orientation of $s$ to be from $Q_1$ to $Q_2$. Since $\tau(g)=-g$, and
$\tau(W')=W'$, we deduce that $\tau(s)=-s$, and $Q_1, Q_2$ are interchanged by $\tau$. It follows that $Q_1$ and $Q_2$ are both contained in
either $g_1$ or $g_2$. Without loss of generality, we can assume that $Q_1,Q_2$ are contained in $g_2$.

\noindent We know that $\tau$ preserves the orientation of $M$, since $\tau$ reverses the orientation of $s$ and $g_2$, we deduce that $s$ meets
both sides of $g_2$. Let $r$ be one of the two subsegments of $g_2$ with endpoints $Q_1,Q_2$, and let $W''$ be the Weierstrass point which is
contained in $r$. Note that $Q_1$ and $Q_2$ must be distinct, otherwise $s=g$, and $\card\{g\cap(g_1\cup g_2)\}= 1$, which is discard by the
hypothesis. Consider the simple closed curve $g'_3$ which is composed by $s$ and $r$. Note that we have $\langle g'_3, g_2 \rangle=\pm 1$, therefore $g'_3\neq 0$ in $H_1(M,\Z)$, and in particular $g'_3$ is an essential, non-separating curve. From the definition, we see that $\tau(g'_3)=-g'_3$. We can move $g'_3$ slightly in its free homotopy class so that the following conditions are satisfied

\begin{itemize}
\item[.] $g'_3\cap s =\{W'\}$,

\item[.] $g'_3 \cap g_2 =\{W''\},$

\item[.] $\tau(g'_3)=-g'_3$.
\end{itemize}

\noindent By construction, we have

\begin{itemize}

\item[.] $g'_3\cap g_1=\vide.$

\item[.] $\card\{g'_3\cap g\}=\card\{\inter(r) \cap g \}+1 \leq \card\{g_2\cap g\} -2+1 =\card\{g_2\cap g\}-1.$

\end{itemize}

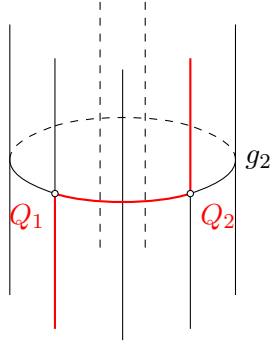
\begin{figure}[htb]
\begin{center}
\begin{tikzpicture}[scale=0.3]

\draw (-5,0) .. controls (-5,-0.5) and (-4.5,-1) .. (-3,-1.5); \draw (5,0) .. controls (5,-0.5) and (4.5,-1) .. (3,-1.5);

\draw[dashed] (-5,0) .. controls (-5,0.5) and (-4.5,1) .. (-3,1.5); \draw[dashed] (5,0) .. controls (5,0.5) and (4.5,1) .. (3,1.5);
\draw[dashed] (-3,1.5) .. controls (-1.5,2) and (1.5,2) .. (3,1.5);

\draw (0,4) -- (0,-8); \draw[dashed] (-1,7) -- (-1,-4) (1,7) -- (1,-4);

\draw[thick, red] (-3,-1.5) -- (-3,-7.5) (3,-1.5) -- (3,4.5); \draw (-5,6) -- (-5,-6) (5,6) -- (5,-6); \draw (-3,-1.5) -- (-3,4.5) (3,-1.5) --
(3,-7.5);

\draw[red, thick] (-3,-1.5) .. controls (-1.5,-2) and (1.5,-2) .. (3,-1.5);

\filldraw[fill=white] (-3,-1.5) circle (4pt) (3,-1.5) circle (4pt);

\draw[red] (-3,-1.5) node[below left] {$Q_1$} (3,-1.5) node[below right] {$Q_2$}; \draw (5,0)  node[right] {$g_2$};
\end{tikzpicture}
\end{center}
\caption{ Case $W' \in g_3$}
\label{fig:ReduceIntNb:Case1}
\end{figure}

Let $g_3$ be the simple closed geodesic in the free homotopy class of $g'_3$. Since $\tau(g'_3)=-g'_3$, we have $\tau(g_3)=-g_3$, as $\tau(g_3)$ is the simple closed geodesic in the free homotopy class of $\tau(g'_3)$.  It follows from Lemma~\ref{lm:simple-arc:sphere} that $g_3$ does not contain $W$, therefore we can conclude that $g_3$ verifies the property $(\P)$.

Let us now show that $g_3$ satisfies the conditions in the conclusion of the lemma. Let $\iota$ denote the geometric intersection number between free homotopy classes of simple closed curves in $M$.  Recall that $\iota(\alpha,\beta)=\card\{\alpha_0,\beta_0\}$, where $\alpha_0$ and $\beta_0$ are thesimple closed geodesics in the free homotopy classes of $\alpha$ and $\beta$ respectively. We have

\begin{itemize}

\item[.] $\card\{g_3\cap g_1\}=\iota(g'_3,g_1) =0$,

\item[.] $\card\{g_3\cap g_2\}=\iota(g'_3,g_2) \leq \card\{g'_3\cap g_2\}=1$,

\item[.] $\card\{g_3 \cap g \}=\iota(g'_3,g) \leq \card\{g'_3\cap g\} < \card\{g_2\cap g\}$.
\end{itemize}

\noindent By construction, we have $\langle g_3,g_2\rangle=\langle g'_3,g_2\rangle=\pm 1$, therefore $g_3\cap g_2 \neq \vide$. We deduce that
$\card\{g_3\cap g_2\}=1$, and the lemma is proven for this case.

\bigskip

\noindent \underline{\bf Case 2: $W' \notin g$.} Cutting $M$ along the curves $g_1,g_2$, we then get $4$-holed sphere $N$ which is equipped with
a hyperbolic metric with geodesic boundary. Let $\hat{g}$ denote the union of geodesic arcs with endpoints in $\partial N$ corresponding to sub-segments of $g$ with endpoints in $g_1\cup g_2$. Let $\hat{s}_1$ be a geodesic arc realizing the distance $\di_N(W',\hat{g})$. Note that $\hat{s}_1$ does not meet $\partial N$. The involution $\tau$ of $M$ induces an involution on $N$, which will be denoted by $\tau_N$. Let $\hat{s}_2$ denote the image of $\hat{s}_1$ by $\tau_N$. Note that $\hat{s}_2$ is also a geodesic arc realizing the distance $\di_N(W',\hat{g})$. 

From the fact that both $\hat{s}_1, \hat{s}_2$ realize the distance in $N$ from $W'$ to $\hat{g}$, we derive that  $W'$ is the unique common point of $\hat{s}_1$ and $\hat{s}_2$, since if any other common point exists, it must be one endpoint of both segments, hence it is fixed by $\tau_N$. But $\tau_N$ has only two fixed points in $N$ corresponding to $W$ and $W'$, and by assumption $W\not\in g$, and we get a contradiction. Since $\tau$ acts like $-\Id$ on the tangent plan at $W'$, $\hat{s}=\hat{s}_1\cup\hat{s}_2$ is in fact a geodesic segment.

Let $s$ be the geodesic arcs in $M$ corresponding to $\hat{s}$, and let $Q_i, \; i=1,2,$ denote the endpoint of $s$. By construction, we have $s\cap (g_1\cup g_2)=\vide$, $W \not\in  s$, and $\tau(s)=-s$. Let $R_1$ be an endpoint of the segment of $g$ (with endpoints in $g_1\cap g_2$) that contains $Q_1$, and let $r_1$ denote the oriented subsegment from $Q_1$ to $R_1$. We note by $R_2$ and $r_2$ the images of $R_1$ and $r_1$ by the involution $\tau$ respectively. The curve $c=(-r_1)*(-s_1)*s_2*r_2$ is then a simple arc joining $R_1$ to $R_2$ verifying $\tau(c)=-c$.

Since  $\tau(R_1)=R_2$, it follows that either $\{R_1,R_2\}\subset g_1$, or $\{R_1,R_2\}\subset g_2$. Without loss of generality, we
can assume that $R_1,R_2$ are contained in $g_2$, then the same argument as in {\bf Case 1} shows that $c$ meets both sides of $g_2$. Here we
have two issues

\begin{itemize}

\item[$\bullet$] $R_1=R_2$: in this case $c$ is actually a simple closed curve which satisfies

\begin{itemize}

\item[a)] $\tau(c)=-c$.

\item[b)] $c\cap g_1=\vide$.

\item[c)] $\card\{c\cap g_2\}=1$. 

\end{itemize}

We can then find in the neighborhood of $c$ a simple closed curve $c'$, freely homotopic to $c$, satisfying a), b), c) such that $c'\cap g'_3=\{R_1\}$. If $\card(g\cap g_2)>1$,  let $g_3$ be $c'$. If $\card(g\cap g_2)=1$, that is $g\cap g_2=\{R_1\}$, then we take $g'_3$ to be the image of $c'$ by the Dehn twist about $g_2$, observe that $g'_3$ satisfies a), b), c), and does not meet $g$ (see Figure~\ref{fig:ReduceIntNb:Case2}). We denote by $g'_3$ the simple closed geodesic in the free homotopy class of $g'_3$.

\begin{figure}[htb]
\begin{center}
\begin{tikzpicture}[scale=0.4]

\draw (-4,0) arc (180:360: 4 and 1.5); \draw[dashed] (4,0) arc (0:180: 4 and 1.5); 

\draw (-4,4) -- (-4,-4) (4,4) -- (4,-4); \draw (0,-1) -- (0,-2); \draw[red, thick] (0,2.5) -- (0,-1) (0,-2) -- (0,-5);

\draw[red, thick] (0,-1) arc (270:360: 4 and 1.5); \draw[red, thick, dashed] (4,0.5) .. controls (4,3) and (-4,1) .. (-4,-0.5); \draw[red,
thick] (-4,-0.5) arc (180:270: 4 and 1.5);

\filldraw[fill=white] (0,-1.5) circle (3pt); \draw (0,-1.5) node[below right] {$R_1$} (4,0) node[right] {$g_2$};

\end{tikzpicture}
\end{center}
\caption{ Case $W' \not\in g$ and $R_1=R_2$.}
\label{fig:ReduceIntNb:Case2}
\end{figure}
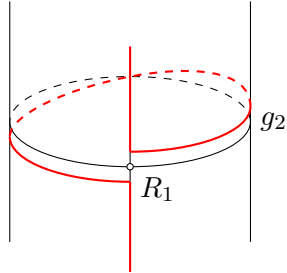

\item[$\bullet$] $R_1\neq R_2$: let $d$ be an arc in $g_2$ with endpoints $R_1,R_2$, then $c\cup d$ is a simple closed curve invariant under $\tau$.
We can find in the neighborhood of $c\cup d$ a simple closed curve $g'_3$ which satisfies a), b), c), and $\card\{g'_3\cap g\} \leq \card\{g_2\cap g\}-2$. We can then take $g_3$ to be the closed geodesic in the free homotopy class of $g'_3$.

\end{itemize}

\noindent In both cases, using the same arguments as in {\bf Case 1}, we see that $g_3$ verifies the required properties. It is also clear from the  arguments above that, if in addition, we have $g\cap g_1=\vide$, then $g_3\cap g_1=\vide$. The proof of the lemma is now complete. \carre

\begin{lemma}\label{lm:simple-arc:sphere}
Let $c$ be a simple closed curve in $M$ such that $\tau(c)=-c$ and $W\not\in c$. Let $c_0$ be the simple closed geodesic in the free homotopy class of $c$, then we also have $\tau(c_0)=-c_0$, and $W\not\in c_0$.
\end{lemma}
\dem Since $\tau(c)=-c$, the image $\bar{c}$ of $c$ under $\rho$ is a simple arc with endpoints in $\{P_1,\dots,P_5\}$ such that $\inter(\bar{c})\cap \{P_0,\dots,P_5\}=\vide$. Assume that $P_1,P_2$ are the endpoints of $\bar{c}$. Let $\tilde{c}$ be a simple closed curve in $\CP^1$ based at $P_0$, which bounds a disc ${\rm D}$ containing $\bar{c}$ such that $\inter({\rm D})\cap\{P_0,\dots,P_5\}=\{P_1,P_2\}$. Let $\tilde{d}$ be a simple arc contained in ${\rm D}$ which joins $P_0$ to an endpoint of $\bar{c}$. We have

\begin{itemize}
\item[.] $\rho^{-1}(\tilde{c})$ is the union of two simple closed curve $c_1,c_2$, freely homotopic to $c$, such that $c_1\cap c_2=\{W\}$.

\item[.] $\rho^{-1}(\tilde{d})$ is a simple closed curve $d$ based at $W$ such that $\tau(d)=-d$. 

\item[.] $\rho^{-1}(\inter({\rm D}))$ is an open annulus bounded by $c_1,c_2$ which contains $d\setminus\{W\}$.
\end{itemize}

\noindent Therefore, in the group $\pi_1(M,W)$, we have $[\tau(c_1)]=[c_2]^{-1}=[d][c_1]^{-1}[d]^{-1}$, $[\tau(d)]=[d]^{-1}$. Clearly, in the homology level, we have $\langle c,d \rangle = \pm 1$. The lemma then follows from the same arguments as in Lemma~\ref{PrfTh1Lm1}.  \carre 

Using Lemma~\ref{Key1Lm}, we can now prove

\begin{proposition}\label{PrfTh1Prop1}
Let $\De=\{a,b_1,b_2,c_1,c_2,d\}$ and $\De'=\{a',b'_1,b'_2,c'_1,c'_2,d'\}$ be two admissible decompositions of $(M,W)$. As usual, let
$(a,b,c,e)$ and $(a',b',c',e')$ be the symplectic bases of $H_1(M,\Z)$ associated to $\De$ and $\De'$ respectively. Assume that $b_1$ is freely
homotopic to $b'_1$, then there exists $\gamma \in \Gamma$ such that

$$\left(%
\begin{array}{c}
  a' \\ b' \\ c' \\ e' \\ \end{array}%
\right)=\gamma \cdot \left(%
\begin{array}{c}
  a \\ b \\ c \\ e \\
\end{array}%
\right).$$

\end{proposition}

\dem Let $b_0,c_0,c'_0$ be the simple closed geodesics in the free homotopy classes of $b_1,c_1,c'_1$ respectively. According to Lemma
\ref{StdDecpLm}, we know that $b_0,c_0,c'_0$ verify Property $(\P)$, and $b_0\cap c_0=b_0\cap c'_0=\vide$.

If $c'_0=\pm c_0$, the proposition follows from Corollary \ref{StdDecpCor}. Hence, we only need to consider the case where $c'_0\neq c_0$ as
subsets of $M$. In this case, since each of the curves $b_0,c_0,c'_0$ contains exactly two Weierstrass points of $M$, and $W \notin b_0\cup
c_0\cup c'_0$, we deduce that $c'_0\cap c_0\neq \vide$. Let $n$ be the number of intersection points of $c_0$ and $c'_0$. The proposition will
be proved by induction.\\

\noindent \underline{\bf Case $n=1$:} Note that in this case, the intersection point between $c_0$ and $c'_0$ must be a Weierstrass point. We
will show that there exist two admissible decompositions $\{\hat{a},\hat{b}_1,\hat{b}_2,\hat{c}_1,\hat{c}_2,\hat{d}\}$, and
$\{\hat{a}',\hat{b}'_1,\hat{b}'_2,\hat{c}'_1, \hat{c}'_1,\hat{c}'_2,\hat{d}'\}$ such that $\hat{b}_i=\hat{b}'_i$, and $\hat{b}_i,\hat{c},\hat{c}'$ are freely homotopic to $b_0$, $c_0$, and $c'_0$ respectively.  We can then use Proposition \ref{PrfTh1Prop0} to conclude. 

\noindent Observe that $\bar{b}_0=\rho(b_0), \bar{c}_0=\rho(c_0), \bar{c}'_0=\rho(c'_0)$ are three simple arcs on $\CP^1$, which satisfy

\begin{itemize}

\item[.] $\rho(W)\notin (\bar{b}_0\cup \bar{c}_0 \cup \bar{c}'_0)$.

\item[.] $\bar{b}_0\cap (\bar{c}_0\cup \bar{c}'_0)=\vide$.

\item[.] $\bar{c}_0$ and $\bar{c}'_0$ have a common endpoint, and  $\card\{\bar{c}_0\cap \bar{c}'_0\}=1$.

\end{itemize}

\begin{figure}[htb]
\begin{center}
\begin{tikzpicture}[scale=0.5]

\draw[red] (0,0) .. controls (-2,3) and (-7,3) .. (-7,0); \draw[red] (0,0) .. controls (-2,-3) and (-7,-3) .. (-7,0);

\draw[red] (0,0) .. controls (3,3) and (8,3) .. (8,0); \draw[red] (0,0) .. controls (3,-3) and (8,-3) .. (8,0);

\draw[red, dashed] (0,0) .. controls (4,-3) and (6,2) .. (9,2); \draw[red, dashed] (12,0) arc (0:90:3 and 2); \draw[red, dashed] (0,0) ..
controls (3,-4) and (12,-4) .. (12,0);

\filldraw[black] (0,0) circle (2pt) (-5,0) circle (2pt) (-3,0) circle (2pt) (4,0) circle (2pt) (7,0) circle (2pt) (10,0) circle (2pt);

\draw (-5,0) -- (-3,0) (4,0) -- (10,0);

\draw (-4,0) node[above] {$\bar{b}_0$} (6,0) node[below] {$\bar{c}_0$} (9,0) node[above] {$\bar{c}'_0$}; \draw (-4,2.2) node[above] {$\tilde{b}$}; \draw (5,2.2) node[above] {$\tilde{c}$}; \draw (12,0) node[right] {$\tilde{c}'$};

\end{tikzpicture}

\end{center}
\caption{Case $b_0=b'_0$, $\card\{c_0\cap c'_0\}=1$.}
\label{fig:bcommon:InterNb:1}
\end{figure}
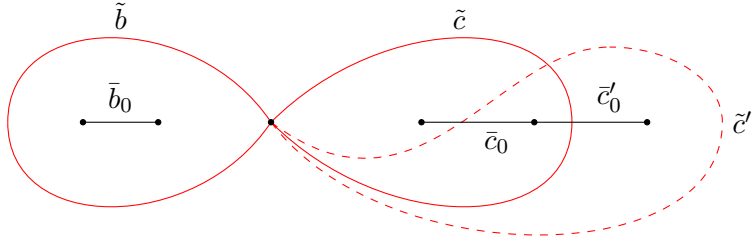

We choose the numbering of $\{P_1,\dots,P_5\}$ so that $P_1$ and $P_2$ are the endpoints of $\bar{b}_0$, $P_3$ and $P_4$ are the endpoints of $\bar{c}_0$, and $P_4$ and $P_5$ are the endpoints of $\bar{c}'_0$. Let $\tilde{b}$ be a simple closed curve in  $\CP^1\sm (\bar{b}_0\cup \bar{c}_0 \cup \bar{c}'_0)$, with base-point $P_0$, which separates $\bar{b}_0$ from $\bar{c}_0\cup \bar{c}'_0$. Let $\mathrm{D}_1$ denote the open disk in $\CP^1$ bounded by $\tilde{b}$ which contains $\bar{b}_0$. The pre-image of $\tilde{b}$ in $M$ consists of two simple closed curves with base point $W$, denoted by $\hat{b}_1$ and $\hat{b}_2$, which bound an open annulus containing $b_0$. The orientation of $b_0$ induces an orientation of $\hat{b}_1$ and $\hat{b}_2$ by free homotopy.  

Let $\tilde{c}$ (resp. $\tilde{c}'$) be a simple closed curve in $\CP^1\setminus \mathrm{D}_1$ with base-point $P_0$ surrounding $\bar{c}_0$ (resp. $\bar{c}'_0$) as shown in Figure~\ref{fig:bcommon:InterNb:1}. Like $\rho^{-1}(\tilde{b})$,  $\rho^{-1}(\tilde{c})$ (resp. $\rho^{-1}(\tilde{c}')$) is the union of two simple closed curves $\hat{c}_1, \hat{c}_2$ (resp. $\hat{c}'_1,\hat{c}'_2$) which bound an embedded open annulus in $M$ containing $c_0$ (resp. $c'_0$). We denote by ${\rm D}_1$ (resp. by ${\rm D}_2$) the open disc in $\CP^1$ bounded by $\tilde{c}$ (resp. by $\tilde{c}'$) and containing $\bar{c}_0$ (resp. $\bar{c}'_0$).

Observe that the family of curves $\{\hat{b}_1, \hat{b}_2,\hat{c}_1,\hat{c}_2\}$ (resp. $\{\hat{b}_1,\hat{b}_2, \hat{c}'_1,\hat{c}'_2\}$) cuts $M$ into three pieces: two annulus and a quadrilateral.  We can then add to these families some simple closed curves to obtain two admissible decompositions for $(M,W)$ as follows: let $\tilde{a}$ be a simple arc in ${\rm D}_1$ joining $P_0$ to and endpoint of $\bar{b}_0$, and let $\tilde{d}$ (resp. $\tilde{d}'$) denote a simple arc in ${\rm D}_2$ (resp. in ${\rm D}'_2$) joining  $P_0$ to an endpoint of $\bar{c}_0$  (resp. of $\bar{c}'_0$).  Set $\hat{a}=\rho^{-1}(\tilde{a}),\hat{d}=\rho^{-1}(\tilde{d})$, and $\hat{d}'=\rho^{-1}(\tilde{d}')$. By choosing appropriate orientations for $\hat{a},\hat{d}$, and $\hat{d}'$, we see that $\hat{\De}=\{\hat{a},\hat{b}_1,\hat{b}_2,\hat{c}_1,\hat{c}_2, \hat{d}\}$, and
$\hat{\De}'=\{\hat{a},\hat{b}_1,\hat{b}_2,\hat{c}'_1,\hat{c}'_2, \hat{d}'\}$ are two admissible decompositions for $(M,W)$. Note that the orientations of $\hat{c}_i$ and $\hat{c}'_i$ are determined by the orientation of $\hat{b}_i$.

Clearly, by construction, we have $\hat{c}_i$ (resp. $\hat{c}'_i$) is freely homotopic to $\pm c_0$ (resp. to $\pm c'_0$). Therefore, by Corollary \ref{StdDecpCor}, the symplectic bases of $H_1(M,\Z)$ associated to $\De$ and $\hat{\De}$ are related by an element of $\Gamma$.  Similarly, the symplectic bases of $H_1(M,\Z)$ associated to $\De'$ and $\tilde{\De}'$ are also related by an element of $\Gamma$. Now, by Proposition \ref{PrfTh1Prop0}, we know that the symplectic bases associated to $\hat{\De}$ and $\hat{\De}'$ are related by an element of $\Gamma$. Hence, the proposition is proven for this case.

\bigskip

\noindent \underline{\bf Case  $n>1$:} By Lemma \ref{Key1Lm}, there exists a simple closed geodesic $c''_0$ verifying Property $(\P)$ such that

\begin{itemize}

\item[.] $\card\{c''_0\cap b_0\}=0$,

\item[.] $\card\{c''_0\cap  c_0\}=1$,

\item[.] $\card\{c'_0 \cap c''_0\} < \card(c'_0\cap c_0)=n$.

\end{itemize}

\noindent Let $\De''=\{a'',b''_1,b''_2,c''_1,c''_2,d''\}$ be the standard decomposition associated to the pair of simple closed geodesics
$(b_0,c''_0)$. The arguments in  Case $n=1$ show that the symplectic bases associated to $\De''$ and $\De$ are related by an element of $\Gamma$.
Now, since $\card\{c'_0\cap c''_0\} < n$, the induction hypothesis implies that the symplectic bases associated to $\De''$ and $\De'$ are also
related by an element of $\Gamma$, and the proposition follows. \carre

\subsection{Proof of Theorem~\ref{th1}}

We can now give the proof of Theorem~\ref{th1}.  Let $b_0,c_0$ be the simple closed geodesics in the free homotopy classes of $b_1$ and $c_1$ respectively. Note that the roles of $b_1$ and $c_1$ are interchanged by the $S$ move, and the orientation of $b_1$ is irrelevant because $-\Id =S^2\in \Gamma$. Let $b'_0$ denote the simple closed geodesic in the free homotopy class of $b'_1$. If $b'_0=\pm b_0$, or $b'_0=\pm c_0$, then Proposition~\ref{PrfTh1Prop1} allows us to conclude immediately.  Assume that $b'_0\neq \pm b_0$ and $b'_0\neq \pm c_0$. Let $n$ be the number of intersection points of $b'_0$ and $b_0\cup c_0$. We  proceed by induction.

\medskip

\noindent \underline{\bf Case $n=1$:} in this case, we can suppose that $b'_0\cap b_0=\vide$, and $\card\{b'_0\cap c_0\}=1$. Let $\De''$ be the standard decomposition associated to the pair $(b_0,b'_0)$. By Proposition~\ref{PrfTh1Prop1}, we know that the symplectic bases of $H_1(M,\Z)$ associated to $\De$ and to $\De''$ are related by  an element of $\Gamma$, and the symplectic bases associated to $\De'$ and $\De''$ are also related by  an element of $\Gamma$, Theorem~\ref{th1} is then proven for this case. 

\medskip

\noindent \underline{\bf Case $n>1$:} by Lemma~\ref{Key1Lm}, there exists a simple closed geodesic $g$ verifying Property $(\P)$ such that

\begin{itemize}

\item[.] $\card\{g\cap b_0\}=0$,

\item[.] $\card\{g\cap c_0\}=1$,

\item[.] $\card\{b'_0\cap g\}< \card\{b'_0\cap c_0\}$.

\end{itemize}

\noindent Let $\De''$ be the standard decomposition associated to  the pair of geodesics $(b_0,g)$. From Proposition~\ref{PrfTh1Prop1}, we know that the symplectic bases of $H_1(M,\Z)$ associated to $\De$ and $\De''$ are related by an element of $\Gamma$. Since $\card\{b'_0\cap (b_0\cup g)\} < \card\{b'_0\cap (b_0\cup c_0)\}$, by the induction hypothesis, the  symplectic bases of $H_1(M,\Z)$ associated to $\De''$ and $\De'$ are related by an element $\Gamma$. Theorem~ref{th1} is then proven. \carre

\section{Proof of Theorem \ref{th0}}\label{sect:proofTh0}

\subsection{The map $\Xi$}

Let $\mathscr{M}$ denote the quotient $\H(2)/\C^*$. We define the map $\Xi$ from $\mathscr{M}$ to $\Gamma\backslash\Hb_2$ as follows: given
a pair $(M,W)$ in $\mathscr{M}$, we associate to $(M,W)$ the $\Gamma$-orbit of the period matrix of the symplectic homology basis associated to an admissible decomposition for $(M,W)$. It follows from Theorem \ref{th1}  that the map $\Xi$ is well defined. We will show that $\Xi$ is a homeomorphism from $\H(2)/\C^*$ and $\Gamma\backslash\mathcal{J}_2$, which implies Theorem \ref{th0}.

\subsection{Injectivity of $\Xi$}

Let $(M,W)$ and $(M',W')$ be two pairs in $\mathscr{M}$. Assume that $M$ and $M'$ are defined by  the equations
$\DS{w^2=\prod_{i=0}^5(z-\lbd_i)}$, and $\DS{w^2=\prod_{i=0}^5(z-\lbd'_0)}$ so that $W$, and $W'$ correspond to $\lbd_0$ and $\lbd'_0$
respectively.  Let $\Pi$ (resp. $\Pi'$) be the period matrix of the symplectic basis associated to an admissible decomposition $\De$  (resp. $\De'$)  for the pair $(M,W)$ (reps. $(M',W'))$.  Assume that there exists an element $\gamma$ of $\Gamma$ such that $\Pi'=\gamma\cdot \Pi$. By Lemma~\ref{lm:Gam:Adm}, we know that there exists an admissible decomposition $\tilde{\De}$ for the pair $(M,W)$ such that the symplectic homology bases of $M$ associated to $\tilde{\De}$ and $\De$ are related by $\gamma$. It follows that the period matrix of the basis associated to $\tilde{\De}$ is equal to $\Pi'$.

\noindent Using an element of $\mathrm{PSL}(2,\C)$ we can assume that $\lbd_1=\lbd'_1=0, \lbd_2=\lbd'_2=1$,  and $\lbd_0=\lbd'_0=\infty$. Then
from Theorem~\ref{ThetThm}, we see that the values of $\lbd_i$ and $\lbd'_i$ ($i=3,4,5$) can be computed by the same theta functions, with the
same period matrix. Therefore we have $\lbd_i=\lbd'_i, \; i=3,4,5$, and it follows that there exists a conformal homeomorphism $\phi: M \ra M'$ such that $\phi(W)=W'$.

\subsection{Surjectivity of $\Xi$}

Let $\tilde{\Pi}$ be a matrix in $\mathcal{J}_2$, we will show that there exists a pair $(M,W)$ in $\mathscr{M}$ such that $\Xi((M,W))=\Gamma\cdot
\tilde{\Pi}$. Since $\tilde{\Pi}\in \mathcal{J}_2$, there exists a Riemann surface of genus two $M$ and a symplectic homology basis whose
period matrix is $\tilde{\Pi}$. Let $W_0$ be a Weierstrass point of $M$, and let $\De$ be an admissible decomposition for the pair $(M,W_0)$. Let $\Pi$ be the period matrix of the symplectic homology basis of $M$ associated to $\De$. By definition, there exists an element $A\in \mathrm{Sp}(4,\Z)$ such that $\tilde{\Pi}=A\cdot \Pi$. According to Lemma~\ref{lm:Sp:Adm}, there exists an admissible decomposition $\tilde{\De}$ for a pair $(M,W)$, where $W$ is also a Weierstrass point of $M$, such that the symplectic homology bases associated to $\De$ and $\tilde{\De}$ is related by $A$.  Consequently, the period matrix of the basis associated to $\tilde{\De}$ is $\tilde{\Pi}$, and by definition $\Gamma\cdot
\tilde{\Pi}=\Xi((M,W))$.


\subsection{Continuity of $\Xi$}

To prove the continuity of $\Xi$ we will consider the inverse map $\Xi^{-1}:\; \Gamma\backslash \mathcal{J}_2 \ra \mathscr{M}$.  Let $\Pi$ be a matrix in $\mathcal{J}_2$,  then $\Pi$ is the period matrix of the symplectic homology basis associated to an admissible decomposition for a pair $(M,W)$ in $\mathscr{M}$.  There exist complex numbers $\{\lbd_0,\lbd_1,\dots,\lbd_5\}$ such that $M$ is the surface defined by the equation $\DS{w^2=\prod_{i=0}^5(z-\lbd_i)}$. We can assume that $W$ is the Weierstrass point corresponding to $\lbd_0$. A neighborhood of $\Pi$ in $\Hb_2$ consists of period matrices of the same symplectic homology basis on Riemann surfaces close to $M$.

\noindent Using $\mathrm{PSL}(2,\C)$, we can assume that $\lbd_0=0$, it follows that $\DS{\omega=\frac{zdz}{w}}$  is a holomorphic $1$-form with
double zero at $W$. Let $\De=\{a,b_1,b_2,c_1,c_2,d\}$ be an admissible decomposition for the pair $(M,W)$.  Let $(a,b,c,e)$ be the symplectic homology basis associated to $\De$, then the map


$$\begin{array}{cccc}
\Phi:& \mathcal{U} & \lra &\C^4\\
     & (M,\omega)  & \longmapsto & (\int_{a}\omega, \int_{b}\omega, \int_{c}\omega, \int_{e}\omega)\\

\end{array}$$

\noindent is a local chart for $\H(2)$ in the neighborhood $\mathcal{U}$ of $(M,\omega)$.  Let $\rho : M \rightarrow \CP^1$ be the two-sheeted branched cover from $M$ onto $\CP^1$. Recall that by construction, $\rho(a)=\bar{a}, \rho(b)=\bar{b}_*, \rho(c)=\bar{c}_*, \rho(e)=\bar{e}$ are simple arcs in $\CP^1$ with endpoints in $\{\lbd_0,\dots,\lbd_5\}$. We have

\begin{eqnarray*}
\int_a\omega & = & 2\int_{\bar{a}}\frac{zdz}{w},\\
\int_b\omega & = & 2\int_{\bar{b}_*}\frac{zdz}{w},\\
\int_c\omega & = & 2\int_{\bar{c}_*}\frac{zdz}{w},\\
\int_e\omega & = & 2\int_{\bar{e}}\frac{zdz}{w}.\\
\end{eqnarray*}
\noindent Clearly, the integrals of $\DS{\frac{zdz}{w}}$ along $\bar{a},\bar{b}_*,\bar{c}_*,\bar{e}$ depend continuously on $(\lbd_1,\dots,\lbd_5)$. Since $\lbd_i$ can be computed from $\Pi$ by some theta functions, we get a continuous map $\Psi$ from a neighborhood of $\Pi$ in $\mathcal{J}_2$ into $\C^4$. Now, in a neighborhood of $\Pi$,  the map $\Xi^{-1}$ is the composition of $\Psi$ and the natural projection $\C^4\sm \{0\} \lra \CP^3$.
It follows immediately that $\Xi^{-1}$ is continuous. Since $\dim_{\C}\mathcal{J}_2=\dim_{\C}\CP^3=3$, we conclude that $\Xi$ is a homeomorphism.

\subsection{$[\mathrm{Sp}(4,\Z):\Gamma]=6$}

We have a natural projection from $\mathscr{M}$ onto the moduli space of closed Riemann surface of genus two $\mathfrak{M}_2$, which is
homeomorphic to $\mathrm{Sp}(4,\Z)\backslash \J_2$, by associating to any element $(M,W)$ of $\mathscr{M}$ the point $M$ in $\mathfrak{M}_2$.
Every Riemann surface of genus two has six Weierstrass points, and the group of automorphisms of a generic one contains exactly two elements,
the identity and the hyperelliptic involution, both fix all the Weierstrass  points. Therefore, the pre-image of a generic point in
$\mathfrak{M}_2$ contains exactly six points. Note that $-\Id_4 \in \Sp(4,\Z)$ acts trivially on $\Hb_2$, and the action of $\Sp(4,\Z)/\{\pm \Id_4\}$ on $\Hb_2$ is effective. Since $\{\pm \Id_4\} \subset \Gamma$, we derive $[\mathrm{Sp}(4,\Z):\Gamma]=6$. The proof Theorem \ref{th0} is now complete. \carre

\bigskip

\begin{appendices}

\section{Existence of parallelogram decompositions}\label{app:sect:ParDecomp:Existence}

\begin{proposition}[Existence of Parallelogram Decompositions]\label{prop:Ext:ParDecomp}
For any  translation surface $(M,\omega)$ in $\H(2)$, there always exists a parallelogram decomposition on $M$.
\end{proposition}

\dem Recall that a {\em simple cylinder} in $(M,\omega)$ is a cylinder such that each boundary component is a saddle connection. We first prove that there always exists a simple cylinder in $M$. By a well known theorem of Masur (see \cite{Ma}, or \cite{MasTab}), there always exists a closed geodesic $\gamma$ in $M$. Let $C_\gamma$ be the (open) cylinder consisting of closed geodesics freely homotopic to $\gamma$. Since $M$ is of genus two, and has only one singular points, we have three cases:

\begin{itemize}

\item[a)] $C_\gamma$ is a simple cylinder.

\item[b)] $M\setminus\overline{C_\gamma}$ is a simple cylinder.

\item[c)] $\overline{C_\gamma}=M$.

\end{itemize}

\begin{figure}[htb]
\begin{center}
\begin{tikzpicture}[scale=0.6]
\filldraw[black] (0,0) circle (2pt) (3,0) circle (2pt) (5,0) circle (2pt) (10,0) circle (2pt); \filldraw (2,4) circle (2pt) (7,4) circle (2pt)
(9,4) circle (2pt) (12,4) circle (2pt);

\draw[thick] (0,0) -- (10,0) -- (12,4) -- (2,4) -- cycle;

\draw (1.5,0) node[above] {$\gamma_1$}; \draw (10.5,4) node[above] {$\gamma_1$}; \draw (4,0) node[above] {$\gamma_2$}; \draw (8,4) node[above]
{$\gamma_2$}; \draw (7.5,0) node[above] {$\gamma_3$}; \draw (4.5,4) node[above] {$\gamma_3$};

\draw (1,2) -- (11,2); \draw (3,2) node[above] {$\gamma$};

\draw[dashed] (7,4) -- (3,0); \draw[dashed] (9,4) -- (5,0);

\end{tikzpicture}
\end{center}
\caption{Surface admitting a decomposition in one cylinder.}
\label{fig:Decomp:1cyl}
\end{figure}
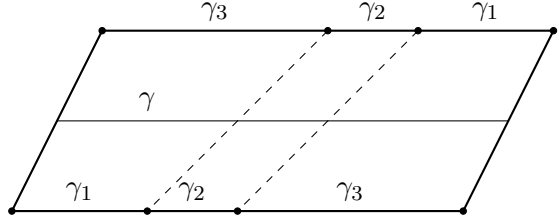

\noindent If $a)$ or $b)$ occurs then we are done. If we are in case $c)$, then the pair $(M,\omega)$ is obtained from  a single parallelogram
as shown in Figure~\ref{fig:Decomp:1cyl}. Therefore, in this situation, one can easily find a simple cylinder in another direction.

\medskip

\rem This argument has been used in the proof of Theorem 7.1 in \cite{McM07}.

\medskip


Let $C_1$ be a simple cylinder in $M$. Let $b_1,b_2$ denote the two saddle connections which bound $C_1$. Cutting off $C_1$ from $M$ along $b_1$
and $b_2$, we get a torus minus two open disks whose boundaries meet at one point, which is the unique singular point $p$ of $M$. Splitting the point $p$ into two points, we obtain then a flat torus with an open disk removed, whose boundary is composed by two
parallel segments of same length. Gluing these two segments together, we obtain a flat torus $M_1$ with a marked simple geodesic arc $b$. 
Let $p_1,p_2$ denote the endpoints of $b$ in $M_1$. Let us show that there always exist two parallel closed geodesic $c_1, c_2$ on $M_1$ such
that $c_1\cap b=\{p_1\},$ and $c_2\cap b=\{p_2\}$.

Choose a direction $\theta\in[0;2\pi]$ which is not parallel to $b$, and consider the geodesic flow $\psi_\theta$ of $M_1$ in this
direction. Let $t_0$ be the smallest value of $t$ such that $t_0>0$ and $\psi^{t_0}_\theta(b)\cap b \neq \vide$. The value $t_0$ must be finite
because otherwise the area of stripe $S_t=\bigcup_{0\leq s\leq t}\psi^s_\theta(b)$ tends to infinity as $t \rightarrow \infty$, which is impossible.

Since $\psi^{t_0}_\theta(b)$ and $b$ are parallel and have the same length, $\psi^{t_0}_\theta(b)$ contains at least one of the endpoints
$p_1,p_2$. Without loss of generality, we can assume that $\psi^{t_0}_\theta(b)$ contains $p_1$. Remark that the stripe
$S_{t_0}=\bigcup_{0\leq s \leq t_0}\psi^s_\theta(a)$ the image of a parallelogram $P$ under an isometric immersion $\varphi$, which is injective in the interior. By assumption, the inverse image of $p_1$ by $\varphi$ contains two points which belong to two opposite sides of $P$. It follows that the image of the segment joining these two points in $\varphi^{-1}(p_1)$ under $\phi$ is a closed geodesic $c_1$ on $M_1$ which verifies $c_1\cap b=\{p_1\}$. Consequently, the closed geodesic $c_2$ passing through $p_2$ and parallel to $c_1$ also verifies $c_2\cap b=\{p_2\}$.

By construction, we can identify $M_1\setminus b$ with the complement of $\overline{C_1}$ in $M$. By this identification, the closed geodesics $c_1$
and $c_2$ become two saddle connections which bound an open cylinder $C_2$ disjoint from $C_1$, and the complement of
$\overline{C_1}\cup\overline{C_2}$ is an open disk isometric to an open parallelogram in $\R^2$.

Let $a$ be a saddle connection contained in $\overline{C_1}$ transverse to the closed geodesics parallel to $b_1$ and $b_2$, let $d$ be a saddle connection in $\overline{C_2}$ transverse to the closed geodesics parallel to  $c_1$ and $c_2$. One can easily check that the family $\{a,b_1,b_2,c_1,c_2,d\}$, with appropriate orientations, is a parallel decomposition of $M$.\carre

It follows from this proposition that, on any surface $(M,\omega)$ in $\H(2)$, there exist infinitely many parallelogram decompositions, since if we have one, we can get infinitely many others by using elementary moves $T,S,R$. It is also possible to show that any two parallelogram decompositions are related by elementary moves.

\section{ Proof of Lemma \ref{SpGenLm}}\label{prfSpGenLm}

For any $g\geq 1$, let $\sigma$ be the permutation of $\{1,2,\dots,2g\}$ that transpose $2i$, and $2i-1$, for $i=1,\dots,g$. The {\em elementary
symplectic matrices} are the matrices

$$E_{ij}=\left\{%
\begin{array}{ll}
    \Id_{2g}+e_{ij}, & \hbox{ if $i=\sigma(j)$;} \\
    \Id_{2g} + e_{ij} -(-1)^{i+j}e_{\sigma(j)\sigma(i)}, & \hbox{ otherwise,} \\
\end{array}%
\right.$$

\noindent where $i\neq j$, and $e_{ij}$ is the matrix whose the $(i,j)$-entry is $1$, and all other entries are $0$. It is a classical fact that
$\mathrm{Sp}(2g,\Z)$ is generated by elementary symplectic matrices (\cite{FarMar}, Chap. 7). For the case $g=2$, we have

\begin{itemize}

\item[$\bullet$] $E_{12}=E^t_{21}=\left(%
\begin{array}{cccc}
  1 & 1 & 0 & 0 \\
  0 & 1 & 0 & 0 \\
  0 & 0 & 1 & 0 \\
  0 & 0 & 0 & 1 \\
\end{array}%
\right), E_{34}=E^t_{43}=\left(%
\begin{array}{cccc}
  1 & 0 & 0 & 0 \\
  0 & 1 & 0 & 0 \\
  0 & 0 & 1 & 1 \\
  0 & 0 & 0 & 1 \\
\end{array}%
\right).$

\item[$\bullet$] $E_{13}=E_{42}^{-1}=\left(%
\begin{array}{cccc}
  1 & 0 & 1 & 0 \\
  0 & 1 & 0 & 0 \\
  0 & 0 & 1 & 0 \\
  0 & -1 & 0 & 1 \\
\end{array}%
\right), \; E_{31}=E_{24}^{-1}=\left(%
\begin{array}{cccc}
  1 & 0 & 0 & 0 \\
  0 & 1 & 0 & -1 \\
  1 & 0 & 1 & 0 \\
  0 & 0 & 0 & 1 \\
\end{array}%
\right).$

\item[$\bullet$] $E_{14}=E_{32}=\left(%
\begin{array}{cccc}
  1 & 0 & 0 & 1 \\
  0 & 1 & 0 & 0 \\
  0 & 1 & 1 & 0 \\
  0 & 0 & 0 & 1 \\
\end{array}%
\right), E_{41}=E_{23}=\left(%
\begin{array}{cccc}
  1 & 0 & 0 & 0 \\
  0 & 1 & 1 & 0 \\
  0 & 0 & 1 & 0 \\
  1 & 0 & 0 & 1 \\
\end{array}%
\right).$

\end{itemize}

\noindent All we need is to verify that $E_{ij}, (i\neq j)$ is contained in the group $\Gamma'$ generated by $\{T,R,S,U\}$. It is clear that $E_{12}, E_{34}, E_{43}$ belong to $\Gamma\subset \Gamma'$. We have

\begin{itemize}

\item[$\bullet$] $E_{21}=U^{-1}T^{-1}U \in \Gamma'$.

\item[$\bullet$] $E_{13}=SU \in \Gamma'$.
\end{itemize}

\noindent Since $\Gamma$ contains $\left(%
\begin{array}{cc}
  \Id_2 & 0 \\
  0 & SL(2,\Z) \\
\end{array}%
\right)$, we have $S_1=\left(%
\begin{array}{cccc}
  1 & 0 & 0 & 0 \\
  0 & 1 & 0& 0 \\
  0 & 0 & 0 & -1 \\
  0 & 0 & 1 & 0 \\
\end{array}%
\right) \in \Gamma$,  therefore $S_2=U^{-1}S_1U=\left(%
\begin{array}{cccc}
  0& -1 & 0 & 0 \\
  1 & 0 & 0 & 0 \\
  0 & 0 & 1 & 0 \\
  0 & 0 & 0 & 1 \\
\end{array}%
\right) \in \Gamma'$. It follows

\begin{itemize}
\item[$\bullet$] $E_{31}=(S_1S_2)E_{13}^{-1}(S_1S_2)^{-1} \in \Gamma'$.

\item[$\bullet$] $E_{14}=S_1E_{13}S_1^{-1} \in \Gamma'$.

\item[$\bullet$] $E_{41}=S_2E_{13}S_2^{-1} \in \Gamma'$.

\end{itemize}

\noindent The lemma is then proven. \carre

\section{ A family of $\Gamma$ right cosets in $\mathrm{Sp}(4,\Z)$}\label{GamRCst}

In this section, we explicit a partition of the group $\mathrm{Sp}(4,\Z)$ into $\Gamma$ right cosets. Recall that, see Lemma \ref{SpGenLm}, the
group $\mathrm{Sp}(4,\Z)$ is generated by $T,S,R$, and $U$. Set

$$ \mathcal{F}=\{\Gamma, U\cdot \Gamma, RU\cdot\Gamma, SRU\cdot\Gamma, URU\cdot \Gamma, USRU\cdot\Gamma\}.$$

\noindent By Lemma \ref{GamPropLm}, we know that the action of $\Gamma$ on $(\Z/2\Z)^4\sm \{0\}$ has two orbits $\mathscr{O}_1$ and
$\mathscr{O}_2$, therefore we have a simple criterion to show that an element of $\mathrm{Sp}(4,\Z)$ {\bf does not} belong to $\Gamma$.
Consequently, it is easy to verify that the elements in the family $\mathcal{F}$ are all distinct.

\noindent We will also determine explicitly the action of $T^{\pm 1}, R^{\pm 1}, S^{\pm 1}, U^{\pm 1}$ on $\mathcal{F}$ by  multiplication from
the left. Note that, since $S^{-1}=-S$ (resp. $ U^{-1}=-U$), and $-\Id_4\in \Gamma$, the actions of $S$ and $S^{-1}$ (resp. $U$ and $U^{-1}$)
are identical. Details of the calculations are lengthly and uninteresting, hence will be omitted. The final result is resumed in the following
table. \\

\begin{center}

\begin{tabular}{|c|c|c|c|c|c|c|}
  \hline
    & $\Gamma$ & $U\cdot\Gamma$ & $RU\cdot\Gamma$ & $SRU\cdot\Gamma$ & $URU\cdot\Gamma$ & $USRU\cdot\Gamma$ \\
  \hline
  $T$ & $\Gamma$ & $U\cdot\Gamma$ & $RU\cdot\Gamma$ & $URU\cdot\Gamma$ & $SRU\cdot\Gamma$ & $USRU\cdot\Gamma$ \\
  \hline
  $R$ & $\Gamma$ & $RU\cdot\Gamma$ & $U\cdot\Gamma$ & $SRU\cdot\Gamma$ & $URU\cdot\Gamma$ & $USRU\cdot\Gamma$ \\
  \hline
  $S$ & $\Gamma$ & $U\cdot\Gamma$ & $SRU\cdot\Gamma$ & $RU\cdot\Gamma$ & $USRU\cdot\Gamma$ & $URU\cdot\Gamma$ \\
  \hline
  $U$ & $U\cdot\Gamma$ & $\Gamma$ & $URU\cdot\Gamma$ & $USRU\cdot\Gamma$ & $RU\cdot\Gamma$ &  $SRU\cdot\Gamma$ \\
  \hline
  $T^{-1}$ & $\Gamma$ & $U\cdot\Gamma$ & $RU\cdot\Gamma$ & $URU\cdot\Gamma$ & $SRU\cdot\Gamma$ & $USRU\cdot\Gamma$ \\
  \hline
  $R^{-1}$ & $\Gamma$ & $RU\cdot\Gamma$ & $U\cdot\Gamma$ & $SRU\cdot\Gamma$ & $URU\cdot\Gamma$ & $USRU\cdot\Gamma$ \\
  \hline
\end{tabular}

\end{center}

%
%
%

\end{appendices}

\bigskip


\begin{thebibliography}{22}
\addcontentsline{toc}{section}{References}

\bibitem{Bus} P. Buser: {\it Geometry and Spectra of Compact Riemann Surfaces}. Progress in Mathematics, Birkhäuser (1992).

\bibitem{EskOk} A. Eskin, A. Okounkov: Asymptotics of number of branched coverings of a torus and volume of moduli spaces of
holomorphic diferentials. Invent. Math., \textbf{145:1}, 59-104 (2001).


\bibitem{FarMar} B. Farb, D. Margalit: {\it A Primer on Mapping Class Groups}.

\bibitem{FarKr} H. Farkas, I. Kra:  {\it Riemann surfaces}, second edition. Graduate Texts in Mathematics, 71, Springer-Verlag, New York,
(1992).

\bibitem{Ha} R. Hain: Finiteness and Torelli spaces, in Problems on Mapping Class Groups and Related Topics, edited by Benson Farb (2006),
Proc. Symp. Pure Math., Amer. Math. Soc. (2006).

\bibitem{HubZo} P. Hubert, H. Masur, T.A. Schmidt, A. Zorich: Problems on billiards, flat surfaces and translation surfaces. In collection:
Problems on Mapping Class Groups and Related Topics, edited by B. Farb, Proc. Symp. Pure Math., Amer. Math. Soc.(2006).

\bibitem{Kon} M. Konsevich: Lyapunov exponents and Hodge theory. ``The mathematical beauty of physics" (Saclay, 1996), (in Honor of C.
Itzykson) 318-332, Adv. Ser. Math. Phys., 24. World Sci. Publishing, River Edge, NJ(1997).

\bibitem{KonZo}M. Konsevich, A. Zorich: Connected components of the moduli spaces of Abelian differentials. Invent. Math., {\bf 153:3}, 631-678
(2003).

\bibitem{Ma} H. Masur: Closed trajectories for quadratic differentials with an application to billiards. Duke Math. J., {\bf 53}, 307-314,
(1986).

\bibitem{MasTab} H. Masur, S. Tabachnikov: Rational billards and flat structures. In: B. Hasselblatt and A. Katok (ed): Handbook of Dynamical
Systems, Vol. 1A, Elsevier Sience B.V., 1015-1089 (2002).


\bibitem{McM03} C.T. McMullen:  Billiards and Teichmüller curves on Hilbert modular surfaces. J. Amer. Math. Soc. {\bf 16}, 857-885 (2003).

\bibitem{McM07} C.T. McMullen: Dynamics of $SL_2(\R)$ over moduli space in genus two. Ann. of Math.(2), 165 No.2, 397-456 (2007).

\bibitem{McM09} C.T. McMullen: Braid groups and Hodge theory ({\it preprint}).

\bibitem{Mes} G. Mess: The Torelli  groups for genus 2 and 3 surfaces. Topology, {\bf 31}, No. 4, 775-790 (1992).










\bibitem{RauFar} H. Rauch, H. Farkas: {\it Theta functions with applications to Riemann surfaces}, The Williams \& Wilkins Company, Baltimore
(1974).


\bibitem{Zor} A. Zorich: Flat surfaces. In collection ``Frontiers in Number Theory, Physics and Geometry", Vol. 1: On random matrices, zeta
functions and dynamical systems, Ecole de Physique des Houches, France, March 9-21 2003, Springer-Verlag (2006).


\end{thebibliography}
\end{document}